\documentclass[a4paper,10pt,reqno]{amsart}

\usepackage[english]{babel}
\usepackage[utf8]{inputenc}
\usepackage{amsmath,amssymb,amsfonts,amsthm,amscd}
\usepackage[mathscr]{eucal}
\usepackage{hyperref}
\usepackage{xcolor}
\usepackage{enumitem}
\setlist[itemize]{topsep=0pt,leftmargin=2em}
\setlist[enumerate]{topsep=0pt,leftmargin=2em}

\newcommand{\R}{\mathbb{R}}
\newcommand{\C}{\mathbb{C}}

\newcommand{\df}{\mathrm{d}}
\newcommand{\Nil}{\mathrm{Nil}_3}
\newcommand{\Iso}{\mathop{\rm Iso}\nolimits}
\DeclareMathOperator{\arcsinh}{arcsinh}

\DeclareMathOperator{\arctanh}{arctanh}
\DeclareMathOperator{\Ric}{Ric}
\DeclareMathOperator{\area}{area}
\DeclareMathOperator{\vol}{vol}
\DeclareMathOperator{\trace}{trace}

\newtheorem{theorem}{Theorem}[section]
\newtheorem{proposition}[theorem]{Proposition}
\newtheorem{corollary}[theorem]{Corollary}
\newtheorem{lemma}[theorem]{Lemma}

\theoremstyle{definition}
	\newtheorem{example}[theorem]{Example}
\theoremstyle{remark}
	\newtheorem{remark}[theorem]{Remark}

\numberwithin{equation}{section}

\title[Total curvature and area growth of entire minimal graphs]{Total curvature and area growth of entire minimal graphs in Heisenberg space}
\date{}

\author{José M. Manzano}
\address{Departamento de Matemáticas, Universidad de Jaén -- Jaén (Spain)}
\email{jmprego@ujaen.es}



\thanks{Research supported by MCIN/AEI/10.13039/501100011033 grant no.\ PID2025-174972NB-I00.}

\begin{document}

\begin{abstract}
We prove that entire minimal graphs in Heisenberg space $\mathrm{Nil}_3(\tau)$ have infinite total curvature. Since their Gauss curvature is negative, this means that their intrinsic area growth is not quadratic, as conjectured by Pérez, Rodríguez and the author~\cite{MPR}. It follows that vertical planes are the only orientable complete stable surfaces with constant mean curvature and quadratic area growth in $\mathrm{Nil}_3(\tau)$, and this classification is extended to surfaces with critical and supercritical constant mean curvature in the homogeneous spaces $\mathbb{E}(\kappa,\tau)$. We also prove that the height of an entire minimal graph in $\mathrm{Nil}_3(\tau)$ grows at most quadratically, as conjectured by Nelli and the author~\cite{ManNel}, whence the extrinsic and cylindrical area growths of such a graph are exactly cubic. Our arguments rely on the conformal Calabi--Lee duality between minimal graphs in $\mathrm{Nil}_3(\tau)$ and spacelike graphs with constant mean curvature $\tau$ in the Lorentz--Minkowski space $\mathbb{L}^3$, for which we also get sharp gradient estimates.
\end{abstract}

\maketitle

\section{Introduction}

Simply connected homogeneous Riemannian 3-manifolds with $4$-dimensional isometry group are classified as a family $\mathbb{E}(\kappa,\tau)$ with two parameters $\kappa,\tau\in\R$ such that $\kappa-4\tau^2\neq 0$. They are characterized by admitting a global unit Killing vector field whose integral curves are the fibers of a Riemannian submersion $\pi:\mathbb{E}(\kappa,\tau)\to\mathbb{M}^2(\kappa)$ with constant bundle curvature $\tau$ over the complete simply connected surface $\mathbb{M}^2(\kappa)$ with constant curvature $\kappa$, see~\cite{Dan,Man}. If $\tau=0$, the family reduces to product spaces $\mathbb{M}^2(\kappa)\times\R=\mathbb{E}(\kappa,0)$, whereas we find the Heisenberg group $\Nil$, the special unitary group $\mathrm{SU}(2)$ and the universal cover of the special linear group $\widetilde{\mathrm{SL}}_2(\mathbb{R})$, all of them endowed with appropriate left invariant metrics, if $\tau\neq 0$.

A cornerstone result in the theory of surfaces with constant mean curvature $H$ ($H$-surfaces in the sequel) in $\mathbb{E}(\kappa,\tau)$-spaces was the solution to the Bernstein problem in $\Nil(\tau)=\mathbb{E}(0,\tau)$ by Fernández and Mira~\cite{FM}. In contrast to the case of Euclidean space $\mathbb{R}^3$, in which all entire minimal graphs are affine planes, in $\Nil(\tau)$ there is a generically $2$-parameter family of entire minimal graphs (i.e., global minimal sections of the Killing submersion $\pi:\Nil(\tau)\to\mathbb{R}^2$) with Abresch--Rosenberg differential~\cite{AR2,Dan2} any prescribed holomorphic quadratic differential $Q$ on a non-compact simply connected Riemann surface $\Sigma$ (with the unique restriction that $Q$ is not identically zero if $\Sigma$ is parabolic). This description is very satisfactory from the point of view of the classification, and it is specially suited to the study of entire minimal graphs from a conformal viewpoint~\cite{DH,Dan2}. However, it does not give much information about their geometry.

The first question in this direction was raised by Pérez, Rodríguez and the author~\cite{MPR}, who proved that an orientable parabolic complete stable $H$-surface immersed in $\Nil(\tau)$ is either a vertical plane (i.e., the preimage of a straight line of $\mathbb{R}^2$ by the Killing submersion) or an entire minimal graph. Since surfaces with quadratic area growth (i.e., the area of their intrinsic metric balls at some point grows at most quadratically with respect to their radius) are parabolic~\cite{CY1} and the area growth of all known examples of entire minimal graphs is cubic, they conjectured that no entire minimal graph in $\Nil(\tau)$ has quadratic area growth (this conjecture was also listed as open problem no.\ 7 in~\cite{DHM}). Later on, the author~\cite{Man19} employed the Calabi--Lee conformal duality between minimal surfaces in $\Nil(\tau)$ and spacelike surfaces with constant mean curvature $\tau$ in the Lorentz--Minkowski space $\mathbb{L}^3$~\cite{Lee,Man19}, along with Cheng and Yau's curvature estimates~\cite{CY2}, to prove that any entire minimal graph $\Sigma\subset\Nil(\tau)$ has negative Gauss curvature. By comparison with the Euclidean plane, it follows that an entire minimal graph in $\Nil(\tau)$ has at least quadratic intrinsic area growth. Li~\cite{Li} showed that a complete surface with negative curvature and quadratic area growth has finite total curvature, and the converse was proved by Hartman~\cite{Hartman}, so that $\Sigma$ has quadratic intrinsic area growth if and only if it has finite total curvature, see~\cite[Thm.~3.4]{Man19}. All in all, the conjecture in~\cite{MPR} can be rephrased in terms of curvature, which is the approach we will follow to prove our first main result.

\begin{theorem}\label{thm:FTC}
There are no entire minimal graphs with finite total curvature or quadratic area growth in $\Nil(\tau)$, with $\tau\neq 0$.
\end{theorem}

The proof of Theorem~\ref{thm:FTC} is based on the analysis of the function $\psi=-\log(\nu)$, where $\nu=\langle N,\xi\rangle>0$ denotes the angle function, being $N$ the upward-pointing unit normal and $\xi$ the unit Killing vector field. The duality reveals that $\psi$ satisfies $0<\Delta\psi\leq -K$, where $K$ is the Gauss curvature of $\Sigma$. If $\Sigma$ had finite total curvature, then it would be conformally equivalent to $\C$, and $\psi$ would become a nonnegative subharmonic function on $\C$ with finite Riesz mass $m$ (see Lemma~\ref{lemma:growth}), so it would grow at most logarithmically. On the other hand, $\nu^{-1}=e^\psi$ is, up to a sign, the angle function of the dual entire $\tau$-graph in $\mathbb{L}^3$, whose Gauss map is a harmonic map into the hyperbolic plane. This fact, together with the growth of $\psi$, enables us to show that its Hopf differential vanishes identically by a Liouville-type argument. However, this Hopf differential is (a constant multiple of) the Abresch--Rosenberg differential of $\Sigma$, which cannot be identically zero if $\Sigma$ is parabolic by Fernández and Mira's classification.

As pointed out in~\cite[Rmk.~2]{MPR}, Theorem~\ref{thm:FTC} leads to the classification of orientable complete stable $H$-surfaces with quadratic area growth in $\Nil(\tau)$, and this result in turn applies in the general framework of $\mathbb{E}(\kappa,\tau)$-spaces, where the same ideas work as long as $4H^2+\kappa\geq 0$ by the Daniel correspondence~\cite{Dan}. Recall that an orientable $H$-surface is said to be (strongly) \emph{stable} if it is a stable critical point of the Jacobi functional $\mathcal{J}=\area-2H\vol$ for all compactly supported normal variations, and that the mean curvature $H$ is called \emph{critical} when $4H^2+\kappa=0$.

\begin{corollary}\label{coro:stable}
Let $\Sigma$ be an orientable complete stable $H$-surface immersed in $\mathbb{E}(\kappa,\tau)$. If $\Sigma$ has quadratic intrinsic area growth, then one of the following holds:
\begin{enumerate}[label=\emph{(\roman*)}]
 \item $\mathbb{E}(\kappa,\tau)=\mathbb{S}^2(\kappa)\times\R$, $H=0$, and $\Sigma=\mathbb{S}^2(\kappa)\times\{t_0\}$ for some $t_0\in\R$;
 \item $4H^2+\kappa=0$ and $\Sigma$ is a vertical cylinder over a complete curve of geodesic curvature $2H$ in $\mathbb{M}^2(\kappa)$.
 \item $4H^2+\kappa<0$ and $\Sigma$ is either a complete vertical graph or a vertical cylinder over a complete curve of geodesic curvature $2H$ in $\mathbb{M}^2(\kappa)$.
\end{enumerate}
In particular, minimal vertical planes are the only orientable complete stable $H$-surfaces with quadratic area growth in $\Nil(\tau)$.
\end{corollary}

Since quadratic area growth implies parabolicity, the results in~\cite{MPR} ensure that a surface $\Sigma$ in the conditions of the statement is either a horizontal slice in $\mathbb{S}^2(\kappa)\times\R$, or a vertical cylinder with critical mean curvature, or a complete vertical multigraph with critical mean curvature (i.e., transversal to the unit Killing vector field). Indeed, Rodríguez and the author~\cite{ManRod} proved that any complete $H$-multigraph in $\mathbb{E}(\kappa,\tau)$ is a graph, which is necessarily entire if the mean curvature is critical, extending the result that complete minimal multigraphs in $\Nil(\tau)$ are entire graphs given by Daniel and Hauswirth~\cite[Thm.~3.1]{DH}. Our contribution is that the case of an entire graph cannot occur under the hypothesis of quadratic area growth, since the Daniel correspondence~\cite{Dan} would give an entire minimal graph in $\Nil(\sqrt{H^2+\tau^2})$ isometric to $\Sigma$, in contradiction with Theorem~\ref{thm:FTC}. 

Ideal Scherk graphs in $\mathbb{H}^2\times\mathbb{R}$ or $\widetilde{\mathrm{SL}}_2(\R)$ are complete stable orientable non-entire $H$-graphs (with subcritical $H$) and quadratic intrinsic area growth, see~\cite[Thm.~4]{ManNel}. We expect that any other complete graph in $\mathbb{H}^2\times\mathbb{R}$ or $\widetilde{\mathrm{SL}}_2(\R)$ has infinite total curvature, which would complete the classification in Corollary~\ref{coro:stable}.

\medskip
The second question we address in this work is related to the other two notions of area growth for a surface $\Sigma$ immersed in $\Nil(\tau)$, proposed by Nelli and the author in~\cite{ManNel}, that is, we can consider the growth of the function $R\mapsto\area(\Sigma\cap A_R)$, where $A_R$ is either the metric ball of $\Sigma$ of radius $R$ (intrinsic area growth), or the metric ball of $\Nil(\tau)$ of radius $R$ (extrinsic area growth), or the solid vertical cylinder of radius $R$, i.e., the preimage of a disk of radius $R$ by the Killing submersion (cylindrical area growth). These concepts do not depend on the point where the sets $A_R$ are centered, and it is easy to see that the intrinsic area growth is slower than the extrinsic one, which is in turn slower than the cylindrical one. If $\Sigma$ is an entire minimal graph in $\Nil(\tau)$, then it was proved in~\cite{ManNel} that the extrinsic (and hence the intrinsic) area growth of $\Sigma$ is at most cubic, whereas its cylindrical area growth is at least cubic and at most quartic. This leads to two subquestions.

On the one hand, since the intrinsic area growth is exactly cubic in all known examples, it was conjectured in~\cite{ManNel} that the intrinsic area growth of \emph{any} entire minimal graph in $\Nil(\tau)$ must be exactly cubic. This question still remains open, since Theorem~\ref{thm:FTC} just says that it is strictly faster than quadratic.

On the other hand, an entire minimal graph $\Sigma_u\subset\Nil(\tau)$ is given by $z=u(x,y)$ in the usual coordinates of $\Nil(\tau)\equiv\R^3$, and it was also conjectured in~\cite{ManNel} that the function $u\in C^\infty(\R^2)$ grows at most quadratically with respect to $r=(x^2+y^2)^{1/2}$, the distance to the origin in the projection to $\R^2$. In \cite[Thm.~6]{ManNel}, it is shown that $|Gu|$ grows at most quadratically and $|u|$ grows at most cubically, where
\begin{equation}\label{eqn:G}
Gu=(u_x+\tau y)\partial_x+(u_y-\tau x)\partial_y
\end{equation}
is the so-called \emph{generalized gradient} of $u$. Nevertheless, it was expected that $|Gu|$ grows at most linearly, and $|u|$ grows at most quadratically, as it happens in all known examples, such as the invariant surfaces of Figueroa, Mercuri and Pedrosa~\cite{FMP} or the graphs of the form $u(x,y)=xf(y)$ given by Daniel~\cite[Ex.~8.4]{Dan2}. It was also shown in~\cite[Cor.~3]{ManNel} that the conjecture on the height growth would imply that $\Sigma_u$ has exactly cubic extrinsic area growth. In the opposite direction, a Collin--Krust type estimate~\cite{CK} was obtained in~\cite[Thm.~7]{ManNel}, which implies that the height of an entire minimal graph grows at least linearly unless it is constant, and this linear lower bound is sharp (as shown by umbrellas or by Cartier's examples~\cite{Cartier}). Our second main result gives a positive answer to the above conjecture.

\begin{theorem}\label{thm:height}
Let $\Sigma_u\subset\Nil(\tau)$ be an entire minimal graph given by $u\in C^\infty(\R^2)$, and let $r=(x^2+y^2)^{1/2}$. Then there exist constants $C_1,C_2>0$ such that
\[|Gu|\leq C_1(1+r)\qquad\text{and}\qquad |u|\leq C_2(1+r^2).\]
\end{theorem}

The constants $C_1$ and $C_2$ are explicit in terms of $\tau$ and the values of $u$ and $Gu$ at one point (see Remarks~\ref{rmk:linear} and~\ref{rmk:height}). Both estimates in Theorem~\ref{thm:height} are sharp, as shown by the entire minimal graph $u(x,y)=\tau xy$. Furthermore, in any entire graph in $\Nil(\tau)$, not necessarily minimal, $|Gu|$ grows at least linearly, see Remark~\ref{rmk:sharp}, so the above estimate says that this growth is exactly linear in the minimal case.

As a consequence of Theorem~\ref{thm:height} and its Lorentzian counterpart (Theorem~\ref{thm:linear}), we are able to determine the exact extrinsic and cylindrical area growths.

\begin{corollary}\label{coro:area}
Entire minimal graphs in $\Nil({}\tau)$ have exactly cubic extrinsic area growth and exactly cubic cylindrical area growth.
\end{corollary}

The estimate for $|Gu|$ in Theorem~\ref{thm:height} is the key of our arguments, since the estimate for $|u|$ follows by integration and the area growths are then essentially deduced from the results in~\cite{ManNel}. It is equivalent, by means of the duality, to the following gradient estimate for entire spacelike graphs with constant mean curvature in $\mathbb{L}^3$, which improves~\cite[Lem.~6]{ManNel} and might be of independent interest.

\begin{corollary}\label{coro:gradient}
Let $\widetilde\Sigma_v\subset\mathbb{L}^3$ be an entire spacelike graph with constant mean curvature $\tau>0$, given by a function $v\in C^\infty(\R^2)$, and let $r=(x^2+y^2)^{1/2}$. Then there exists a constant $C_3>0$ such that
\[|\nabla v|^2\leq 1-\frac{C_3}{1+r^2}.\]
\end{corollary}

In other words, the tangent planes of an entire spacelike $\tau$-graph cannot become lightlike arbitrarily fast as $r\to\infty$. The hyperbolic cylinder $v(x,y)=\frac{1}{2\tau}(1+4\tau^2y^2)^{1/2}$ shows that this estimate is sharp. 

Our approach to these results consists of two steps. First, we use the estimates of Cheng and Yau to show that the vertical sections of $\widetilde\Sigma_v$ are spacelike curves with curvature at most $2\tau$ in a Lorentzian plane, and a comparison with a hyperbola yields the desired gradient estimate provided that $\widetilde\Sigma_v$ has a horizontal tangent plane at some point. Note that in~\cite{ManNel}, the comparison was carried out with straight lines (using convexity), but this is not sharp. Second, we get rid of this assumption by applying a boost of $\mathbb{L}^3$ (i.e., a hyperbolic rotation), which takes a given tangent plane into a horizontal one, see Section~\ref{sec:deformations}. Although boosts are isometries of $\mathbb{L}^3$, they do not preserve the vertical Killing vector field, and they induce a 2-parameter deformation of entire minimal graphs in $\Nil(\tau)$ via the duality. This deformation preserves the conformal structure and the Abresch--Rosenberg differential (see also the appendix of~\cite{Man19}), and changes the metric by a bounded conformal factor that can be controlled in terms of hyperbolic distances between timelike directions.

The paper is organized as follows. In Section~\ref{sec:preliminaries}, we will recall the conformal duality in coordinates and the curvature estimates for entire graphs, we will also introduce the function $\psi$, revisit the harmonic Gauss maps of dual surfaces, and analyze the behavior of the duality under the deformation by boosts. Section~\ref{sec:FTC} will be devoted to proving Theorem~\ref{thm:FTC}, whereas Theorem~\ref{thm:height} and Corollaries~\ref{coro:area} and~\ref{coro:gradient} will be proved in Section~\ref{sec:growth}.

\medskip
\noindent\textbf{AI statement.} This research is the conclusion of a series of papers developed by the author and his collaborators~\cite{MPR,ManRod,ManNel,Man}, containing their original ideas. However, the subharmonicity of the function $\psi$ and the subsequent estimates in Lemma~\ref{lemma:growth} were suggested by an AI agent. Its argument has been simplified, rewritten carefully and corrected at some points by the author, who takes full responsibility for all the results and calculations.

\section{Preliminaries}\label{sec:preliminaries}

\subsection{The conformal duality}\label{sec:duality}
The Heisenberg space $\Nil(\tau)$ and the Lorentz--Minkowski space $\mathbb{L}^3$ will be modeled as $\R^3$ endowed with the metrics
\begin{equation}\label{eqn:ds}
\begin{aligned}
 \df s^2_{\Nil(\tau)}&=\df x^2+\df y^2+(\df z+\tau(y\,\df x-x\,\df y))^2,\\
 \df s^2_{\mathbb{L}^3}&=\df x^2+\df y^2-\df z^2,
\end{aligned}
\end{equation}
respectively. We can define global orthonormal frames $\{E_1,E_2,E_3\}$ in $\Nil(\tau)$ and $\{\widetilde E_1,\widetilde E_2,\widetilde E_3\}$ in $\mathbb{L}^3$, where $\widetilde E_3$ is timelike, as
\begin{equation}\label{eqn:frame}
\begin{aligned}
 E_1&=\partial_x-\tau y\partial_z,& E_2&=\partial_y+\tau x\partial_z,&E_3&=\partial_z,\\
 \widetilde E_1&=\partial_x,&\widetilde E_2&=\partial_y,& \widetilde E_3&=\partial_z.
\end{aligned}\end{equation}
The Killing submersion in both models is nothing but $\pi(x,y,z)=(x,y)$, and the vertical Killing vector field is $\partial_z$. These frames will be assumed positively oriented, which determines a cross product $\times$ in each of the spaces by means of the relation $\langle U\times V,W\rangle=\det(U,V,W)$ for all vector fields $U,V,W$, where the determinant is computed in coordinates with respect to the above positive frames. The Levi-Civita connection $\overline\nabla$ of $\Nil(\tau)$ is given by
\begin{equation}\label{eqn:levi-civita}
\begin{aligned}
\overline\nabla_{E_1}E_1&=0,&
\overline\nabla_{E_1}E_2&=\tau E_3,&
\overline\nabla_{E_1}E_3&=-\tau E_2,\\
\overline\nabla_{E_2}E_1&=-\tau E_3,&
\overline\nabla_{E_2}E_2&=0,&
\overline\nabla_{E_2}E_3&=\tau E_1,\\
\overline\nabla_{E_3}E_1&=-\tau E_2,&
\overline\nabla_{E_3}E_2&=\tau E_1,&
\overline\nabla_{E_3}E_3&=0,
\end{aligned}
\end{equation}
whereas the Levi-Civita connection of $\mathbb{L}^3$ coincides with the usual directional derivative of $\R^3$, which will be denoted by $\mathrm{D}$.

Given a Riemann surface $\Sigma$, a conformal immersion $X:\Sigma\to\Nil(\tau)$ is called \emph{nowhere vertical} if it is transversal to $E_3$, in which case we will consider the upward-pointing unit normal $N$ to the immersion such that the \emph{angle function} $\nu=\langle N,E_3\rangle$ is positive, as well as the shape operator $AU=-\overline\nabla_UN$. Likewise, if $\widetilde X:\Sigma\to\mathbb{L}^3$ is a spacelike conformal immersion, then $\widetilde N$ will denote its future-directed unit normal, i.e., the one whose angle function $\widetilde\nu=\langle\widetilde N,\widetilde E_3\rangle$ is negative, and $\widetilde AU=-\mathrm{D}_U\widetilde N$ will stand for its shape operator. We will assume that the orientation of $\Sigma$ is compatible with these choices of normals, so that the $\frac{\pi}{2}$-rotation $J$ in the tangent bundle is given by $JU=N\times U$ (resp.\ $JU=\widetilde N\times U$). Recall that the mean curvature of $\widetilde X$ is defined as $-\frac{1}{2}\trace(\widetilde A)$ in the Lorentzian setting (see~\cite{Lopez}), in contrast with the mean curvature $\frac{1}{2}\trace(A)$ of $X$ (e.g., the upper sheet of the hyperboloid $x^2+y^2-z^2=-\tau^{-2}$ with Lorentzian radius $\tau^{-1}$ has mean curvature $\tau$).

\medskip
If $\Sigma$ is simply connected, there is a duality between
\begin{enumerate}[label=(\alph*)]
 \item nowhere vertical conformal minimal immersions $X:\Sigma\to\Nil(\tau)$;
 \item spacelike conformal $\tau$-immersions $\widetilde X:\Sigma\to\mathbb{L}^3$.
\end{enumerate}
These dual immersions $X$ and $\widetilde X$ are determined up to a vertical translation, and satisfy $\pi\circ X=\pi\circ\widetilde X$, see~\cite{Lee,Man19}. Since the duality restricts to a duality between graphs over simply connected domains, we will recall how it works in coordinates following~\cite{Man19}. Given an open domain $\Omega\subset\R^2$ and a function $u\in C^\infty(\Omega)$, the vertical graph of $u$ is the surface
\[\Sigma_u=\{(x,y,u(x,y)):(x,y)\in\Omega\}\subset\Nil(\tau).\]
Likewise, any $v\in C^\infty(\Omega)$ satisfying $v_x^2+v_y^2<1$ spans a spacelike graph
\[\widetilde\Sigma_v=\{(x,y,v(x,y)):(x,y)\in\Omega\}\subset\mathbb{L}^3.\]
These graphs are called \emph{entire} if $\Omega=\R^2$. Let us define the quantities
\begin{equation}\label{eqn:alpha-beta}
\begin{aligned}
\alpha&=u_x+\tau y,&\beta&=u_y-\tau x,&\omega&=(1+\alpha^2+\beta^2)^{1/2},\\
\widetilde\alpha&=v_x,&\widetilde\beta&=v_y,&\widetilde\omega&=(1-\widetilde\alpha^2-\widetilde\beta^2)^{1/2}.
\end{aligned}
\end{equation}
Note that $\alpha$ and $\beta$ are the components of the generalized gradient $Gu$ in~\eqref{eqn:G}, so that $\omega^2=1+|Gu|^2$. The graph $\Sigma_u$ is minimal in $\Nil(\tau)$ if and only if
\begin{equation}\label{eqn:minimal}
\partial_x\Big(\frac{\alpha}{\omega}\Big)+\partial_y\Big(\frac{\beta}{\omega}\Big)=0\quad\Longleftrightarrow\quad (1+\beta^2)u_{xx}-2\alpha\beta u_{xy}+(1+\alpha^2)u_{yy}=0,
\end{equation}
whereas the graph $\widetilde\Sigma_v$ has constant mean curvature $\tau$ in $\mathbb{L}^3$ if and only if
\begin{equation}\label{eqn:tau}
\partial_x\Big(\frac{\widetilde\alpha}{\widetilde\omega}\Big)+\partial_y\Big(\frac{\widetilde\beta}{\widetilde\omega}\Big)=2\tau\quad\Longleftrightarrow\quad (1-\widetilde\beta^2)v_{xx}+2\widetilde\alpha\widetilde\beta v_{xy}+(1-\widetilde\alpha^2)v_{yy}=2\tau\widetilde\omega^3.
\end{equation}
In this notation, $\Sigma_u$ and $\widetilde\Sigma_v$ are dual graphs if and only if they satisfy the so-called \emph{twin relations}~\cite{Lee}, namely the equivalent identities
\begin{equation}\label{eqn:twin}
  (\widetilde\alpha,\widetilde\beta)=\Big(\frac{-\beta}{\omega},\frac{\alpha}{\omega}\Big)\quad\Longleftrightarrow\quad
  (\alpha,\beta)=\Big(\frac{\widetilde\beta}{\widetilde\omega},\frac{-\widetilde\alpha}{\widetilde\omega}\Big),
\end{equation}
in which case we also have $\widetilde\omega=\frac{1}{\omega}$. The twin relations give explicit first-order \textsc{pde} systems that allow us to obtain $v$ in terms of $u$ and vice versa:
\begin{align}
v_x&=\tfrac{-(u_y-\tau x)}{\sqrt{1+(u_x+\tau y)^2+(u_y-\tau x)^2}},&v_y&=\tfrac{u_x+\tau y}{\sqrt{1+(u_x+\tau y)^2+(u_y-\tau x)^2}},\label{eqn:twin-v}\\
u_x&=\tfrac{v_y}{\sqrt{1-v_x^2-v_y^2}}-\tau y,&u_y&=\tfrac{-v_x}{\sqrt{1-v_x^2-v_y^2}}+\tau x.\label{eqn:twin-u}
\end{align}
Note that the integrability condition $(v_x)_y=(v_y)_x$ of the system~\eqref{eqn:twin-v} is precisely the minimal surface equation~\eqref{eqn:minimal}, and the integrability condition $(u_x)_y=(u_y)_x$ of the system~\eqref{eqn:twin-u} is nothing but the mean curvature $\tau$ equation~\eqref{eqn:tau}. Therefore, if $\Omega$ is simply connected, then any minimal graph $\Sigma_u$ over $\Omega$ has a dual graph $\widetilde\Sigma_v$ with mean curvature $\tau$ over $\Omega$, and vice versa, each of them being determined by the other one up to an additive constant.

Consider the global frames $\{e_1,e_2\}$ in $\Sigma_u$ and $\{\widetilde e_1,\widetilde e_2\}$ in $\widetilde\Sigma_v$ defined as
\begin{equation}\label{eqn:e1e2}
\begin{aligned}
 e_1&=E_1+\alpha E_3,& \widetilde e_1&=\widetilde E_1+\widetilde\alpha\widetilde E_3,\\
 e_2&=E_2+\beta E_3,& \widetilde e_2&=\widetilde E_2+\widetilde\beta\widetilde E_3,
\end{aligned}\end{equation}
which are the coordinate frames associated with the parametrizations $(x,y)\mapsto(x,y,u(x,y))$ and $(x,y)\mapsto(x,y,v(x,y))$, respectively. The first fundamental forms of $\Sigma_u$ and $\widetilde\Sigma_v$ in these frames are given by
 \begin{equation}\label{eqn:1ff}
 \begin{aligned}
  \mathrm{I}&\equiv\left(\begin{matrix}1+\alpha^2&\alpha\beta\\\alpha\beta&1+\beta^2\end{matrix}\right),&
  \widetilde{\mathrm{I}}&\equiv\left(\begin{matrix}1-\widetilde\alpha^2&-\widetilde\alpha\widetilde\beta\\-\widetilde\alpha\widetilde\beta&1-\widetilde\beta^2\end{matrix}\right),
  \end{aligned}\end{equation}
respectively, so~\eqref{eqn:twin} implies that $\widetilde{\mathrm{I}}=\omega^{-2}\,\mathrm{I}$, whence
\begin{equation}\label{eqn:phi}
\Phi:\Sigma_u\to\widetilde\Sigma_v,\qquad \Phi(x,y,u(x,y))=(x,y,v(x,y)),
\end{equation}
is a global conformal diffeomorphism. In terms of conformal immersions $X:\Sigma\to\Sigma_u$ and $\widetilde X:\Sigma\to\widetilde\Sigma_v$, the condition $\pi\circ X=\pi\circ\widetilde X$ yields $\Phi=\widetilde X\circ X^{-1}$, so the metrics induced by dual immersions are related by $X^*\df s^2_{\Nil(\tau)}=\omega^2\widetilde X^*\df s^2_{\mathbb{L}^3}$. From~\eqref{eqn:e1e2}, we also deduce that the unit normal vector fields to $\Sigma_u$ and $\widetilde\Sigma_v$ are given by
 \begin{equation}\label{eqn:N}
 \begin{aligned}
   N&=-\tfrac{\alpha}{\omega}E_1-\tfrac\beta\omega E_2+\tfrac1\omega E_3,&\widetilde N&=\tfrac{\widetilde\alpha}{\widetilde\omega}\widetilde E_1+\tfrac{\widetilde\beta}{\widetilde\omega}\widetilde E_2+\tfrac{1}{\widetilde\omega} \widetilde E_3.
  \end{aligned}\end{equation}
Hence the twin relations~\eqref{eqn:twin} imply that the angle functions
\begin{equation}\label{eqn:nu}
\nu=\langle N,E_3\rangle=\tfrac{1}{\omega}\qquad\text{and}\qquad\widetilde\nu=\langle \widetilde N,\widetilde E_3\rangle=\tfrac{-1}{\widetilde\omega}
\end{equation}
satisfy $\nu\widetilde\nu=-1$. Note that $\nu$ and $\widetilde\nu$ take values in $(0,1]$ and $(-\infty,-1]$, respectively.

In the sequel, we will denote by $\nabla$, $\Delta$, $\df A$, and $K$ (resp.\ $\widetilde\nabla$, $\widetilde\Delta$, $\widetilde{\df A}$, and $\widetilde K$) the gradient, the Laplacian, the area element, and the Gauss curvature of $\Sigma_u$ (resp.\ $\widetilde\Sigma_v$). Since $\Phi$ induces a conformal factor $\nu^2$ between $\Sigma_u$ and $\widetilde\Sigma_v$, we get that
\begin{equation}\label{eqn:conformal-operators}
\Delta f=\nu^2\,\widetilde\Delta f,\qquad |\nabla f|^2=\nu^2|\widetilde\nabla f|^2,\qquad \df A=\nu^{-2}\,\widetilde{\df A},
\end{equation}
as well as the following relation between Gauss curvatures:
\begin{equation}\label{eqn:conformal-K}
K=\nu^2\widetilde K+\Delta\log(\nu).
\end{equation}
The functions and tensors defined on $\Sigma_u$ and $\widetilde\Sigma_v$ will be identified hereafter through the diffeomorphism $\Phi$ given by~\eqref{eqn:phi}, so the above identities make sense.

\begin{remark}
The twin relations~\eqref{eqn:twin} and the systems~\eqref{eqn:twin-v} and~\eqref{eqn:twin-u} are those in the proof of~\cite[Thm.~2.3]{Man19}, which are compatible with~\eqref{eqn:minimal} and~\eqref{eqn:tau}. We point out that the signs of $\widetilde\alpha$ and $\widetilde\beta$ contain a mistake when the twin relations are introduced in~\cite[\S2.1]{Man19}, and so does the sign of $\widetilde\nu$, which is negative as shown by~\eqref{eqn:nu}. This misprint does not affect the rest of results in~\cite{Man19}.
\end{remark}

\subsection{Curvature estimates}\label{sec:estimates}
Let $\Sigma_u\subset\Nil(\tau)$ and $\widetilde\Sigma_v\subset\mathbb{L}^3$ be dual graphs. Their second fundamental forms are defined as $\sigma(U,V)=\langle\overline\nabla_UV,N\rangle=\langle AU,V\rangle$ and $\widetilde\sigma(U,V)=\langle\mathrm{D}_UV,\widetilde N\rangle=\langle\widetilde AU,V\rangle$, respectively. In the Lorentzian case, we easily get from~\eqref{eqn:e1e2} and~\eqref{eqn:N} that
\begin{equation}\label{eqn:2ff}
\widetilde\sigma(\widetilde e_1,\widetilde e_1)=\frac{-v_{xx}}{\widetilde\omega},\qquad
\widetilde\sigma(\widetilde e_1,\widetilde e_2)=\frac{-v_{xy}}{\widetilde\omega},\qquad
\widetilde\sigma(\widetilde e_2,\widetilde e_2)=\frac{-v_{yy}}{\widetilde\omega}.
\end{equation}
The Gauss equations of $\Sigma_u$ and $\widetilde\Sigma_v$ read $K=\det(A)+\tau^2-4\tau^2\nu^2$ and $\widetilde K=-\det(\widetilde A)$, respectively, see~\cite{Dan,Lopez}. Since $\trace(A)=0$ and $\trace(\widetilde A)=-2\tau$, they can be rewritten as
\begin{equation}\label{eqn:gauss}
|A|^2=-2K+2\tau^2-8\tau^2\nu^2,\qquad|\widetilde A|^2=4\tau^2+2\widetilde K.
\end{equation}
In particular, if we denote by $\widetilde\kappa_1$ and $\widetilde\kappa_2$ the eigenvalues of $\widetilde A$ at some point, then we get the relations $\widetilde\kappa_1+\widetilde\kappa_2=-2\tau$ and $\widetilde\kappa_1\widetilde\kappa_2=-\widetilde K$. We will need the following estimates for entire graphs by Cheng and Yau, see also~\cite{CT,Treibergs}.

\begin{proposition}[\cite{CY2}]\label{prop:CY}
Let $\widetilde\Sigma_v\subset\mathbb{L}^3$ be an entire spacelike $\tau$-graph, with $\tau>0$.
\begin{enumerate}[label=\emph{(\alph*)}]
	\item The induced metric on $\widetilde\Sigma_v$ is complete (and the converse also holds true: any complete spacelike $\tau$-surface immersed in $\mathbb{L}^3$ is an entire graph).
	\item The Gauss curvature of $\widetilde\Sigma_v$ satisfies $-\tau^2\leq\widetilde K\leq 0$ or, equivalently, the shape operator of $\widetilde\Sigma_v$ satisfies $2\tau^2\leq|\widetilde A|^2\leq 4\tau^2$.
\end{enumerate}
In particular, $-2\tau\leq \widetilde\kappa_1,\widetilde\kappa_2\leq 0$ at any point of $\widetilde\Sigma_v$, and hence
\begin{equation}\label{eqn:CY-convex}
0\leq-\widetilde\sigma(V,V)\leq 2\tau\,\langle V,V\rangle
\end{equation}
for all tangent vectors $V$ to $\widetilde\Sigma_v$. It follows from~\eqref{eqn:2ff} that $v$ is a convex function.
\end{proposition}

Observe that the upper bound for $\widetilde K$ in item (b) is the non-trivial one, since the lower bound follows from the inequality $\widetilde\kappa_1\widetilde\kappa_2\leq\frac{1}{4}(\widetilde\kappa_1+\widetilde\kappa_2)^2=\tau^2$, which holds true for any spacelike $\tau$-surface. The last assertion of the statement is not quoted in~\cite{CY2}, but note that $\widetilde\kappa_1$ and $\widetilde\kappa_2$ cannot have opposite signs since $\widetilde\kappa_1\widetilde\kappa_2=-\widetilde K\geq 0$, so they are both non-positive because their sum is $-2\tau<0$.

The angle functions $\nu$ and $\widetilde\nu$ are the normal components of the Killing vector fields $E_3$ and $\widetilde E_3$ along the surfaces $\Sigma_u$ and $\widetilde\Sigma_v$. Therefore, they lie in the kernel of the stability operators $L$ of $\Sigma_u$ and $\widetilde L$ of $\widetilde\Sigma_v$, respectively, which can be expressed by means of~\eqref{eqn:gauss} as
\begin{equation}\label{eqn:jacobi}
\begin{aligned}
L&=\Delta+|A|^2+\Ric(N)=\Delta-2K-4\tau^2\nu^2,\\
\widetilde L&=\widetilde\Delta-|\widetilde A|^2=\widetilde\Delta-2\widetilde K-4\tau^2,
\end{aligned}
\end{equation}
where we have also used that $\Ric(N)=-2\tau^2+4\tau^2\nu^2$ in $\Nil(\tau)$. We will write this information in terms of the function
\begin{equation}\label{eqn:psi}
\psi=-\log(\nu)=\log(\omega),
\end{equation}
which takes nonnegative values and will play a fundamental role in our arguments.

\begin{lemma}\label{lemma:psi}
Let $\Sigma_u\subset\Nil(\tau)$ and $\widetilde\Sigma_v\subset\mathbb{L}^3$ be dual graphs. The function $\psi$ given by~\eqref{eqn:psi} satisfies the following identities:
\begin{enumerate}[label=(\alph*)]
	\item $\Delta\psi=\nu^2\widetilde K-K$,
	\item $\Delta\psi=-2K-4\tau^2\nu^2+|\nabla\psi|^2$,
	\item $\Delta\psi=\nu^2|\widetilde A|^2-|\nabla\psi|^2$.
\end{enumerate}
In particular, the Gauss curvatures of $\Sigma_u$ and $\widetilde\Sigma_v$ are related by
\begin{equation}\label{eqn:bestialidad}
-K=\nu^2(4\tau^2+\widetilde K)-|\nabla\psi|^2.
\end{equation}
\end{lemma}

\begin{proof}
Item (a) is just a restatement of~\eqref{eqn:conformal-K}. Since $\nu=e^{-\psi}$ lies in the kernel of the operator $L$ given by~\eqref{eqn:jacobi}, we get that $\Delta\nu=2K\nu+4\tau^2\nu^3$, so item (b) follows from the identity
\[\Delta\psi=-\Delta\log(\nu)=-\tfrac{1}{\nu}\Delta\nu+\tfrac{1}{\nu^2}|\nabla\nu|^2=-2K-4\tau^2\nu^2+|\nabla\psi|^2.\]
Likewise, since $\omega=-\widetilde\nu=e^\psi$ satisfies $\widetilde L\omega=0$, we find that $\widetilde\Delta\omega=|\widetilde A|^2\omega$, whence
\[\widetilde\Delta\psi=\widetilde\Delta\log(\omega)=\tfrac{1}{\omega}\widetilde\Delta\omega-\tfrac{1}{\omega^2}|\widetilde\nabla\omega|^2=|\widetilde A|^2-|\widetilde\nabla\psi|^2.\]
Item (c) follows from multiplying this equation by $\nu^2$ and taking~\eqref{eqn:conformal-operators} into account. Finally, we get~\eqref{eqn:bestialidad} by subtracting (b) from (a). Note that any two of the items imply the third one in view of the second identity in~\eqref{eqn:gauss}.
\end{proof}

Equation~\eqref{eqn:bestialidad} was essentially employed in~\cite{Man19} to prove that entire minimal graphs in $\Nil(\tau)$ have negative Gauss curvature. We will now give a more direct argument, which is based on the following estimate for the gradient of $\psi$.

\begin{lemma}\label{lemma:gradient}
Let $\Sigma_u\subset\Nil(\tau)$ and $\widetilde\Sigma_v\subset\mathbb{L}^3$ be dual graphs, and let $\widetilde\kappa_1$ and $\widetilde\kappa_2$ be the eigenvalues of the shape operator $\widetilde A$ of $\widetilde\Sigma_v$ at some point. Then, at that point,
\begin{equation}\label{eqn:gradient-psi}
|\nabla\psi|^2\leq\nu^2(1-\nu^2)\max\{\widetilde\kappa_1^2,\widetilde\kappa_2^2\}.
\end{equation}
\end{lemma}

\begin{proof}
Let us decompose $\widetilde E_3=\widetilde T-\widetilde\nu\widetilde N=\widetilde T+\omega\widetilde N$, where $\widetilde T$ is tangent to $\widetilde\Sigma_v$. By taking squared norms in $\mathbb{L}^3$, we get that $|\widetilde T|^2=\omega^2-1=\nu^{-2}-1$. Since $\widetilde E_3$ is a constant vector field in $\mathbb{L}^3$, any tangent vector field $V$ to $\widetilde\Sigma_v$ satisfies
\[V(\omega)=-V\bigl(\langle\widetilde N,\widetilde E_3\rangle\bigr)=-\langle\mathrm{D}_V\widetilde N,\widetilde E_3\rangle=\langle\widetilde AV,\widetilde T\rangle=\langle V,\widetilde A\widetilde T\rangle,\]
which means that $\widetilde\nabla\omega=\widetilde A\widetilde T$. From~\eqref{eqn:conformal-operators} and~\eqref{eqn:psi}, we deduce that
\[|\nabla\psi|^2=\nu^2|\widetilde\nabla\psi|^2=\nu^4|\widetilde\nabla\omega|^2=\nu^4|\widetilde A\widetilde T|^2\leq\nu^4\max\{\widetilde\kappa_1^2,\widetilde\kappa_2^2\}\,|\widetilde T|^2,\]
and the statement follows.
\end{proof}

\begin{proposition}\label{prop:K}
Let $\Sigma_u\subset\Nil(\tau)$ be an entire minimal graph, $\tau>0$. If $K$ and $\nu$ are the Gauss curvature and the angle function of $\Sigma_u$, and $\psi=-\log(\nu)$, then
\begin{equation}\label{eqn:K}
\tfrac{1}{2}\bigl(4\tau^2\nu^2-|\nabla\psi|^2\bigr)\leq -K\leq 4\tau^2\nu^2-|\nabla\psi|^2.
\end{equation}
In particular, $K\leq -3\tau^2\nu^4<0$.
\end{proposition}

\begin{proof}
Let $\widetilde\Sigma_v\subset\mathbb{L}^3$ be the dual entire graph, and let $\widetilde\kappa_1$ and $\widetilde\kappa_2$ be the eigenvalues of its shape operator at some point, so that $-2\tau\leq \widetilde\kappa_1,\widetilde\kappa_2\leq 0$ by Proposition~\ref{prop:CY}. The second inequality in~\eqref{eqn:K} follows from~\eqref{eqn:bestialidad} and the estimate $\widetilde K\leq 0$. As for the first one, observe that it is equivalent to $\Delta\psi\geq 0$ by item (b) of Lemma~\ref{lemma:psi}. Since $|\widetilde A|^2\geq\max\{\widetilde\kappa_1^2,\widetilde\kappa_2^2\}$, we can combine Lemmas~\ref{lemma:psi} and~\ref{lemma:gradient} to reach
\begin{equation}\label{prop:K:eqn1}
\begin{aligned}
	\Delta\psi&=\nu^2|\widetilde A|^2-|\nabla\psi|^2=\nu^2(\widetilde\kappa_1^2+\widetilde\kappa_2^2)-|\nabla\psi|^2\\
	&\geq\nu^2(\widetilde\kappa_1^2+\widetilde\kappa_2^2)-\nu^2(1-\nu^2)\max\{\widetilde\kappa_1^2,\widetilde\kappa_2^2\}\geq\nu^4(\widetilde\kappa_1^2+\widetilde\kappa_2^2)\geq 0.
\end{aligned}
\end{equation}
Finally, item (a) of Lemma~\ref{lemma:psi} along with~\eqref{prop:K:eqn1} yield the estimate
\[\begin{aligned}
-K=\Delta\psi+\nu^2\widetilde\kappa_1\widetilde\kappa_2&\geq\nu^4(\widetilde\kappa_1^2+\widetilde\kappa_2^2)+\nu^2\widetilde\kappa_1\widetilde\kappa_2\\
&\geq\nu^4(\widetilde\kappa_1^2+\widetilde\kappa_2^2+\widetilde\kappa_1\widetilde\kappa_2)=\nu^4(4\tau^2-\widetilde\kappa_1\widetilde\kappa_2)\geq 3\tau^2\nu^4,
\end{aligned}\]
where we have used that $\nu^2\geq\nu^4$ (since $0<\nu\leq 1$) and $0\leq \widetilde\kappa_1\widetilde\kappa_2\leq\tau^2$, as well as the identity $\widetilde\kappa_1+\widetilde\kappa_2=-2\tau$.
\end{proof}

\begin{remark}\label{rmk:K}
These estimates are not enough to decide whether or not the Gauss curvature of an entire minimal graph is integrable. They imply that the total curvature is finite provided that $\nu^2$ is integrable, but this is never the case by~\cite[Cor.~4]{ManNel}, where a result of Espinar~\cite{Esp} was sharpened to show that vertical cylinders are the only orientable complete stable $H$-surfaces in $\mathbb{E}(\kappa,\tau)$, with $\tau\neq 0$ and $4H^2+\kappa\geq 0$, whose angle function is square integrable. They also imply that $\nu^4$ is integrable if the total curvature is finite, but $\nu^4$ can be integrable or not, as shown by the examples below. This is the reason why the gradient term $|\nabla\psi|^2$ in~\eqref{eqn:K} should not be discarded. Actually, the proof of Theorem~\ref{thm:FTC} will rely on the fact that this term and the Gauss curvature can be handled together by means of the identities in Lemma~\ref{lemma:psi}.
\end{remark}

\begin{example}[Umbrellas]\label{ex:umbrella}
The entire minimal graph $\Sigma_u\subset\Nil(\tau)$ with $u(x,y)=0$ is rotationally invariant about the $z$-axis. From~\eqref{eqn:alpha-beta} and~\eqref{eqn:twin-v}, we get that $\alpha=\tau y$ and $\beta=-\tau x$, so that $\omega=(1+\tau^2r^2)^{1/2}$ with $r=(x^2+y^2)^{1/2}$, and the dual surface $\widetilde\Sigma_v\subset\mathbb{L}^3$ is the hyperboloid given by $v(x,y)=\frac{1}{\tau}(1+\tau^2r^2)^{1/2}$. Therefore, $\widetilde\Sigma_v$ is totally umbilical, i.e., $\widetilde\kappa_1=\widetilde\kappa_2=-\tau$ and $\widetilde K=-\tau^2$. The angle function and the Gauss curvature of $\Sigma_u$ are given by
\[\nu=\frac{1}{\sqrt{1+\tau^2r^2}},\qquad K=-\frac{3\tau^2+2\tau^4r^2}{(1+\tau^2r^2)^2}=-2\tau^2\nu^2-\tau^2\nu^4,\]
see~\cite[Rmk.~3.5]{Man19}. Since the metric of $\Sigma_u$ reads $\df r^2+r^2(1+\tau^2r^2)\df\theta^2$ in polar coordinates and $\psi=\frac{1}{2}\log(1+\tau^2r^2)$ only depends on $r$, we can easily compute
\[|\nabla\psi|^2=\psi_r^2=\tau^2\nu^2(1-\nu^2),\qquad\Delta\psi=\tau^2\nu^2+\tau^2\nu^4,\]
where the last formula follows from Lemma~\ref{lemma:psi}. Hence the inequalities in~\eqref{eqn:K} are strict in this case, and $K$ is not integrable in spite of the fact that $\nu^4$ is integrable.
\end{example}

\begin{example}[Invariant surfaces]\label{ex:FMP}
Complete minimal surfaces in $\Nil(\tau)$ invariant under a $1$-parameter group of translations were classified by Figueroa, Mercuri and Pedrosa~\cite{FMP} as the family of entire graphs $\Sigma_{u_\theta}$, $\theta\in\R$, where
\begin{equation}\label{eqn:FMP}
u_\theta(x,y)=\tau xy+\frac{\sinh(\theta)}{4\tau}\left(2\tau y\sqrt{1+4\tau^2y^2}+\arcsinh(2\tau y)\right).
\end{equation}
In this case, $\alpha=2\tau y$ and $\beta=\sinh(\theta)(1+4\tau^2y^2)^{1/2}$, whence $\omega=\cosh(\theta)(1+4\tau^2y^2)^{1/2}$, and the dual surface is the entire graph $\widetilde\Sigma_{v_\theta}\subset\mathbb{L}^3$ given by
\begin{equation}\label{eqn:FMP-dual}
v_\theta(x,y)=\frac{1}{2\tau\cosh(\theta)}\sqrt{1+4\tau^2y^2}-\tanh(\theta)x.
\end{equation}
Hence $\widetilde\Sigma_{v_\theta}$ is a hyperbolic cylinder, which is flat and has $\widetilde\kappa_1=0$ and $\widetilde\kappa_2=-2\tau$. The angle function and the Gauss curvature of $\Sigma_{u_\theta}$ are given by
\[\nu=\frac{1}{\cosh(\theta)\sqrt{1+4\tau^2y^2}},\qquad K=\frac{-4\tau^2}{\cosh^2(\theta)(1+4\tau^2y^2)^2}=-4\tau^2\cosh^2(\theta)\nu^4,\]
see~\cite[Rmk.~3.5]{Man19}. Since $\psi=\log(\cosh(\theta))+\frac{1}{2}\log(1+4\tau^2y^2)$ only depends on $y$, it follows from~\eqref{eqn:1ff} that $|\nabla\psi|^2=\frac{1+\alpha^2}{\omega^2}\psi_y^2$, which leads to
\[|\nabla\psi|^2=4\tau^2\nu^2\bigl(1-\cosh^2(\theta)\nu^2\bigr),\qquad\Delta\psi=4\tau^2\cosh^2(\theta)\nu^4,\]
where we have used Lemma~\ref{lemma:psi} again. Therefore, the second inequality in~\eqref{eqn:K} becomes an equality for these surfaces, and neither $K$ nor $\nu^4$ are integrable.
\end{example}

\subsection{Harmonic Gauss maps}\label{sec:gauss}
Let $X:\Sigma\to\Nil(\tau)$ and $\widetilde X:\Sigma\to\mathbb{L}^3$ be dual conformal immersions. The unit normal $\widetilde N=\sum_{k=1}^3\widetilde N_k\widetilde E_k$ of $\widetilde X$ is just the classical Gauss map into the hyperbolic plane $\mathbb{H}^2$, understood as the unit sphere of $\mathbb{L}^3$. The unit normal of $X$ is written likewise as $N=\sum_{k=1}^3N_kE_k$, with angle function $N_3=\nu>0$, and can be identified with the map $(N_1,N_2,N_3):\Sigma\to\mathbb{S}^2_+\subset\mathbb{R}^3$. Since the frame $\{E_1,E_2,E_3\}$ is left-invariant for the Lie group structure of $\Nil(\tau)$, the map $N$ is usually known as the left-invariant Gauss map of $X$, see~\cite{MeePer,DFM} (Meeks, Mira, Pérez and Ros~\cite{MMPR} have also extended this idea to $H$-surfaces in metric Lie groups, where harmonicity is no longer true in general). Here, we will give direct new proofs of some facts we shall need in our arguments via the duality.

\begin{lemma}\label{lemma:gauss}
The Gauss maps of dual conformal immersions $X:\Sigma\to\Nil(\tau)$ and $\widetilde X:\Sigma\to\mathbb{L}^3$ are related by
\begin{equation}\label{eqn:gauss-maps}
(\widetilde N_1,\widetilde N_2,\widetilde N_3)=(\tfrac{N_2}{N_3},\tfrac{-N_1}{N_3},\tfrac{1}{N_3}).
\end{equation}
\end{lemma}

\begin{proof}
We can assume that $X$ and $\widetilde X$ parametrize dual graphs $\Sigma_u$ and $\widetilde\Sigma_v$, respectively. From~\eqref{eqn:N} and the twin relations~\eqref{eqn:twin}, we deduce that
\begin{equation}\label{lemma:gauss:eqn1}
(N_1,N_2,N_3)=\bigl(\tfrac{-\alpha}{\omega},\tfrac{-\beta}{\omega},\tfrac{1}{\omega}\bigr),\qquad (\widetilde N_1,\widetilde N_2,\widetilde N_3)=(\tfrac{\widetilde\alpha}{\widetilde\omega},\tfrac{\widetilde\beta}{\widetilde\omega},\tfrac{1}{\widetilde\omega})=(-\beta,\alpha,\omega),
\end{equation}
and we are done.
\end{proof}

It is well known that the Gauss map of a spacelike surface with constant mean curvature in $\mathbb{L}^3$ is a harmonic map into $\mathbb{H}^2$~\cite[Thm.~10]{Milnor} (see also~\cite{CT}). Essentially, the fact that $\widetilde X$ has constant mean curvature implies that $\widetilde\Delta\widetilde N=|\widetilde A|^2\widetilde N$ as an $\mathbb{L}^3$-valued map, so that the tangent component of $\widetilde\Delta\widetilde N$ to $\mathbb{H}^2$ vanishes identically (the third coordinate of this identity says that $\widetilde\nu$ lies in the kernel of the stability operator $\widetilde L$ given by~\eqref{eqn:jacobi}). The hyperbolic metric in $\mathbb{S}^2_+\subset\mathbb{R}^3$ reads $\frac{1}{z^2}(\df x^2+\df y^2+\df z^2)$ in the usual coordinates, and it is isometric to $\mathbb{H}^2$ by means of the diffeomorphism $(x,y,z)\mapsto(\frac{y}{z},\frac{-x}{z},\frac{1}{z})$. Therefore, Lemma~\ref{lemma:gauss} yields the harmonicity of the Gauss map $N$ for the hyperbolic metric of $\mathbb{S}^2_+$. This result was proved by Daniel~\cite{Dan2}, and we have just provided a shorter proof via the duality. 

Let $z=s+it$ be a conformal parameter on $\Sigma$, so that the metrics induced by $X$ and $\widetilde X$ can be written as $\mu^2|\df z|^2$ and $\widetilde\mu^2|\df z|^2$, respectively, for some positive functions $\mu$ and $\widetilde\mu$ such that $\widetilde\mu=\nu\mu$. By complexifying the tangent bundle, the conformality of $X$ and $\widetilde X$ can be expressed as
\begin{equation}\label{eqn:conformality}
\langle X_z,X_z\rangle=\langle\widetilde X_z,\widetilde X_z\rangle=0,\qquad \langle X_z,X_{\bar z}\rangle=\tfrac{\mu^2}{2},\qquad\langle\widetilde X_z,\widetilde X_{\bar z}\rangle=\tfrac{\widetilde\mu^2}{2},
\end{equation}
in terms of the usual operators $\partial_z=\frac{1}{2}(\partial_s-i\partial_t)$ and $\partial_{\bar z}=\frac{1}{2}(\partial_s+i\partial_t)$. Consider
\begin{equation}\label{eqn:Q}
Q=i\langle \overline\nabla_{X_z}X_z,N\rangle+2\tau\langle X_z,E_3\rangle^2,\qquad \widetilde Q=\langle\widetilde X_{zz},\widetilde N\rangle,
\end{equation}
so that $Q\,\df z^2$ is the Abresch--Rosenberg differential of $X$ and $\widetilde Q\,\df z^2$ is the Hopf differential of $\widetilde X$, both of which are holomorphic~\cite{Dan2}. They have been normalized so that $4Q=Q_{i,2\tau}(\partial_s,\partial_s)$ and $4\widetilde Q=\widetilde Q_{1,0}(\partial_s,\partial_s)$ in the notation of~\cite[\S2.2]{Man19}, where they are shown to satisfy $Q=\widetilde Q$. This fact was known previously, see~\cite{Dan2,DFM}, but here we are adopting a non-standard normalization so that $Q=\widetilde Q$.

Since the fact that $\widetilde X$ has constant mean curvature $\tau$ can be expressed as $\widetilde X_{z\bar z}=\frac{\tau\widetilde\mu^2}{2}\widetilde N$, a standard computation using~\eqref{eqn:conformality} and~\eqref{eqn:Q} shows that
\begin{equation}\label{eqn:Nz}
\widetilde N_z=\tau\widetilde X_z-\tfrac{2}{\widetilde\mu^2}\widetilde Q\widetilde X_{\bar z},
\end{equation}
which easily implies that $\langle \widetilde N_z,\widetilde N_z\rangle=-2\tau\widetilde Q$, that is,
\begin{equation}\label{eqn:hopf}
\begin{aligned}
	((\widetilde N_1)_z)^2+((\widetilde N_2)_z)^2-((\widetilde N_3)_z)^2&=-2\tau\widetilde Q,\\
	((N_1)_z)^2+((N_2)_z)^2+((N_3)_z)^2&=-2\tau N_3^2 Q,
\end{aligned}
\end{equation}
where the second identity follows from the first one by the substitution~\eqref{eqn:gauss-maps}. The two dual equations in~\eqref{eqn:hopf} give another short proof of the well-known fact that $Q$ and $\widetilde Q$ are the Hopf differentials associated with the harmonic maps $N$ and $\widetilde N$, respectively, except for the factor $-2\tau$, see~\cite{FM,Dan2,DFM}. From~\eqref{eqn:Nz}, we also get that
\begin{equation}\label{eqn:energy}
|\widetilde A|^2=\frac{4}{\widetilde\mu^2}\langle\widetilde A\widetilde X_z,\widetilde A\widetilde X_{\bar z}\rangle=\frac{4}{\widetilde\mu^2}\langle\widetilde N_z,\widetilde N_{\bar z}\rangle=2\tau^2+\frac{8|\widetilde Q|^2}{\widetilde\mu^4}.
\end{equation}

\subsection{Behavior under conformal deformations}\label{sec:deformations}
Let us denote by $\Iso_\xi(\Nil(\tau))$ and $\Iso_\xi^+(\mathbb{L}^3)$ the subgroups of $\Iso(\Nil(\tau))$ and $\Iso(\mathbb{L}^3)$ of direct Killing isometries (i.e., those leaving the unit Killing vector field $\xi=E_3=\widetilde E_3$ invariant, and note that all isometries of $\Nil(\tau)$ are direct). For each $F\in\Iso_\xi(\Nil(\tau))$ and $\widetilde F\in\Iso_\xi^+(\mathbb{L}^3)$, there are unique direct isometries $\rho(F),\widetilde\rho(\widetilde F)\in\Iso^+(\R^2)$ such that $\rho(F)\circ\pi=\pi\circ F$ and $\widetilde\rho(\widetilde F)\circ\pi=\pi\circ\widetilde F$. Therefore, $\rho$ and $\widetilde\rho$ define surjective group morphisms with kernel the subgroups of \emph{vertical translations} (i.e., the isometries spanned by $\xi$), so there is an isomorphism
from $\Iso_\xi(\Nil(\tau))/\ker(\rho)$ to $\Iso_\xi^+(\mathbb{L}^3)/\ker(\widetilde\rho)$. This isomorphism can be extended to an isomorphism
\begin{equation}\label{eqn:R}
R:\frac{\Iso(\Nil(\tau))}{\ker(\rho)}\longrightarrow\frac{\Iso_\xi(\mathbb{L}^3)}{\ker(\widetilde\rho)},
\end{equation}
where $\Iso_\xi(\mathbb{L}^3)$ stands for the group of all isometries of $\mathbb{L}^3$ leaving $\xi$ invariant, but possibly reversing the global orientation.

\begin{proposition}[{\cite{Man19}}]\label{prop:Killing-isometries}
If $X:\Sigma\to\Nil(\tau)$ and $\widetilde X:\Sigma\to\mathbb{L}^3$ are dual conformal immersions, then $F\circ X$ and $R(F)\circ\widetilde X$ are also dual for all $F\in\Iso(\Nil(\tau))$.
\end{proposition}

This means that all isometries of $\Nil(\tau)$ are well understood by the duality, but this is still not the full picture, since the group $\Iso_\uparrow(\mathbb{L}^3)$ of isometries of $\mathbb{L}^3$ preserving the time-orientation has dimension $6$ and its elements also preserve the condition of being an entire spacelike $\tau$-graph by Proposition~\ref{prop:CY}.(a). Given an entire minimal graph $\Sigma_u\subset\Nil(\tau)$ with dual graph $\widetilde\Sigma_v\subset\mathbb{L}^3$ and an isometry $S\in\Iso_\uparrow(\mathbb{L}^3)$, the dual graph of $S(\widetilde\Sigma_v)$ is another entire minimal graph in $\Nil(\tau)$ which is not congruent to $\Sigma_u$ in general, but has the same Abresch--Rosenberg differential. This gives rise to a 2-parameter deformation that can be understood as hyperbolic and parabolic rotations of $\mathbb{L}^3$ fixing an ideal boundary point of $\mathbb{H}^2$. These are hidden symmetries of the family of entire minimal graphs in $\Nil(\tau)$ that are only revealed by the duality, and it is difficult to deal with them in coordinates (e.g., see Lee's explanation that Daniel's examples of the form $z=x f(y)$ have the hidden symmetry via the duality in~\cite[Ex.~4]{Lee}). 

The rest of this section will be devoted to analyzing how these deformations affect the geometry of an entire minimal graph. To this end, it suffices to consider hyperbolic rotations, also known as \emph{boosts}.

Given $\eta_1,\eta_2\in\mathbb{H}^2$, let $d_{\mathbb{H}^2}(\eta_1,\eta_2)$ be the hyperbolic distance from $\eta_1$ to $\eta_2$, which satisfies the formula $\cosh d_{\mathbb{H}^2}(\eta_1,\eta_2)=-\langle \eta_1,\eta_2\rangle$ (it plays the role of a Lorentzian angle between unit timelike vectors). Any $\eta\in\mathbb{H}^2$ can be written as
\begin{equation}\label{eq:delta-theta}
\eta=(\sinh(\delta)\cos(\vartheta),\sinh(\delta)\sin(\vartheta),\cosh(\delta)),
\end{equation}
where $\delta=d_{\mathbb{H}^2}(\eta,(0,0,1))\geq 0$ and $\vartheta\in\R$. The linear map $B_\eta:\mathbb{L}^3\to\mathbb{L}^3$ given by
\begin{equation}\label{eqn:boost}
\begin{aligned}
B_\eta(\cos(\vartheta),\sin(\vartheta),0)&=(\cosh(\delta)\cos(\vartheta),\cosh(\delta)\sin(\vartheta),-\sinh(\delta)),\\
B_\eta(-\sin(\vartheta),\cos(\vartheta),0)&=(-\sin(\vartheta),\cos(\vartheta),0),\\
B_\eta(0,0,1)&=(-\sinh(\delta)\cos(\vartheta),-\sinh(\delta)\sin(\vartheta),\cosh(\delta)),
\end{aligned}
\end{equation}
is an orientation-preserving isometry in $\Iso_{\uparrow}(\mathbb{L}^3)$. Note that $B_\eta(\eta)=\widetilde E_3$.

\begin{remark}\label{rmk:boost}
Note that any $S\in\Iso_\uparrow(\mathbb{L}^3)$ whose linear part takes $\eta$ to $\widetilde E_3$ can be written as $S=S_1\circ B_\eta$ for some $S_1\in\Iso_\xi(\mathbb{L}^3)$. Therefore, the boosts $B_\eta$ span the $2$-parameter deformation of entire minimal graphs that preserves the Abresch--Rosenberg differential (up to ambient isometries), and this deformation is naturally parametrized by $\eta\in\mathbb{H}^2$. Note that $\eta$ represents the timelike direction that becomes vertical in the dual entire spacelike $\tau$-graph, which is a convex surface. This gives a neat explanation of most results in the appendix of~\cite{Man19}. 
\end{remark}

\begin{proposition}\label{prop:deformation}
Let $X:\Sigma\to\Sigma_u\subset\Nil(\tau)$ and $\widetilde X:\Sigma\to\widetilde\Sigma_v\subset\mathbb{L}^3$ be dual conformal immersions spanning entire graphs, and let $\widetilde N:\Sigma\to\mathbb{H}^2$ be the Gauss map of $\widetilde X$. Given a unit vector $\eta\in\mathbb{H}^2$, the following assertions hold true:
\begin{enumerate}[label=(\alph*)]
	\item The immersion $\widetilde X_\eta=B_\eta\circ\widetilde X$ is a conformal immersion of the entire spacelike $\tau$-graph spanned by $-\langle\widetilde X,\eta\rangle$. Its Gauss map is $\widetilde N_\eta=B_\eta\circ\widetilde N$ and its angle function is $\widetilde\nu_\eta=\langle\widetilde N,\eta\rangle=-\cosh d_{\mathbb{H}^2}(\widetilde N,\eta)$.
	\item The dual immersion $X_\eta:\Sigma\to\Nil(\tau)$ of $\widetilde X_\eta$ defines an entire minimal graph in $\Nil(\tau)$. The metrics induced by $X$ and $X_\eta$ in $\Sigma$ are conformal with
	\[X_\eta^*\df s^2_{\Nil(\tau)}=\frac{\langle\widetilde N,\eta\rangle^2}{\langle\widetilde N,\widetilde E_3\rangle^2}X^*\df s^2_{\Nil(\tau)},\]
	and the conformal factor lies between $e^{-2\delta}$ and $e^{2\delta}$, where $\delta=d_{\mathbb{H}^2}(\eta,\widetilde E_3)$.
\end{enumerate}
\end{proposition}

\begin{proof}
Since $B_\eta$ is an isometry of $\mathbb{L}^3$ that preserves the time-orientation, we have that $\widetilde X_\eta$ is a conformal spacelike $\tau$-immersion with future-directed unit normal $B_\eta\circ\widetilde N$, which induces the same metric as $\widetilde X$ on $\Sigma$. This metric is complete, so the image of $\widetilde X_\eta$ is an entire graph by Proposition~\ref{prop:CY}.(a). The height function of $\widetilde X_\eta$ is the third coordinate, which can be computed as
\[-\langle \widetilde X_\eta,\widetilde E_3\rangle=-\langle B_\eta\circ\widetilde X,B_\eta(\eta)\rangle=-\langle \widetilde X,\eta\rangle.\] 
The angle function is in turn the negative third coordinate of the normal:
\[\widetilde\nu_\eta=\langle B_\eta\circ\widetilde N,\widetilde E_3\rangle=\langle B_\eta\circ\widetilde N,B_\eta(\eta)\rangle=\langle\widetilde N,\eta\rangle=-\cosh d_{\mathbb{H}^2}(\widetilde N,\eta),\]
and we have proved (a). As for item (b), we have that $X_\eta$ parametrizes an entire minimal graph because $\widetilde X_\eta$ is in turn entire by item (a). Moreover, we can relate the metrics induced by these dual immersions as
\[X_\eta^*\df s^2_{\Nil(\tau)}=\frac{1}{\nu_\eta^2}\widetilde X_\eta^*\df s^2_{\mathbb{L}^3}=\frac{1}{\nu_\eta^2}\widetilde X^*\df s^2_{\mathbb{L}^3}=\frac{\nu^2}{\nu_\eta^2}X^*\df s^2_{\Nil(\tau)},\]
where $\nu_\eta^{-1}=-\widetilde\nu_\eta=-\langle\widetilde N,\eta\rangle$ and $\nu^{-1}=-\widetilde\nu=-\langle\widetilde N,\widetilde E_3\rangle$ by~\eqref{eqn:nu} and item (a). 

To get the bound for the conformal factor, write $\eta$ as in~\eqref{eq:delta-theta} and let $w=(w_1,w_2,w_3)\in\mathbb{L}^3$ be future-timelike. The Cauchy--Schwarz inequality in $\R^2$ gives
\begin{equation}\label{prop:deformation:eqn1}
\begin{aligned}
	\bigl|\langle w,\eta\rangle+w_3\cosh(\delta)\bigr|&=\sinh(\delta)\bigl|w_1\cos(\vartheta)+w_2\sin(\vartheta)\bigr|\\
	&\leq\sinh(\delta)(w_1^2+w_2^2)^{1/2}\leq\sinh(\delta)w_3.
\end{aligned}
\end{equation}
Since $\cosh(\delta)\pm\sinh(\delta)=e^{\pm\delta}$, Equation~\eqref{prop:deformation:eqn1} amounts to $e^{-\delta}w_3\leq-\langle w,\eta\rangle\leq e^{\delta}w_3$, and it suffices to set $w=\widetilde N$, in which case $w_3=\widetilde N_3=-\langle\widetilde N,\widetilde E_3\rangle$.
\end{proof}

\begin{example}\label{ex:boost}
Consider the entire minimal graph $\Sigma_{u_0}$ given by $u_0(x,y)=\tau xy$ as in Example~\ref{ex:FMP}, whose dual graph $\widetilde\Sigma_{v_0}$ is the hyperbolic cylinder $v_0(x,y)=\frac{1}{2\tau}(1+4\tau^2y^2)^{1/2}$, whose Gauss map $\widetilde N=(0,2\tau y,(1+4\tau^2y^2)^{1/2})$ can be computed by means of~\eqref{lemma:gauss:eqn1}. If $\eta=(\sinh(\theta),0,\cosh(\theta))$ for some $\theta\in\R$, then~\eqref{eqn:boost} yields
\[B_\eta(x,y,v_0(x,y))=\bigl(x\cosh(\theta)-v_0(x,y)\sinh(\theta),\,y,\,-x\sinh(\theta)+v_0(x,y)\cosh(\theta)\bigr),\]
so it is easy to check that $B_\eta(\widetilde\Sigma_{v_0})$ is the hyperbolic cylinder $\widetilde\Sigma_{v_\theta}$ given by~\eqref{eqn:FMP-dual}, i.e., the dual of the entire minimal graph $\Sigma_{u_\theta}$ given by~\eqref{eqn:FMP}. Since $\langle\widetilde N,\eta\rangle=-\cosh(\theta)(1+4\tau^2y^2)^{1/2}$ and $\langle\widetilde N,\widetilde E_3\rangle=-(1+4\tau^2y^2)^{1/2}$, the conformal factor in Proposition~\ref{prop:deformation}.(b) is the constant $\cosh^2(\theta)$, so we recover the fact that the surfaces $\Sigma_{u_\theta}$ are intrinsically homothetic (see the proof of~\cite[Prop.~2]{ManNel}). The hyperbolic cylinder $\widetilde\Sigma_{v_0}$ is invariant under the boosts $B_\eta$ with $\eta=(0,\sinh(\theta),\cosh(\theta))$, so the deformation of $\Sigma_{u_0}$ only depends on one parameter, see also~\cite[Thm.~A.3]{Man19}.
\end{example}

\section{Infinite total curvature}\label{sec:FTC}

Let $\Sigma_u\subset\Nil(\tau)$ be an entire minimal graph, and assume that $\tau>0$ throughout this section. The proof that $\Sigma_u$ has infinite total curvature will rely on the properties of the function $\psi=-\log(\nu)$ introduced in Section~\ref{sec:estimates}. We will begin by collecting the estimates for the Laplacian and the gradient of $\psi$ that can be deduced by the arguments of Proposition~\ref{prop:K}.

\begin{lemma}\label{lemma:subharmonic}
Let $\Sigma_u\subset\Nil(\tau)$ be an entire minimal graph with Gauss curvature $K$ and angle function $\nu$. Then the function $\psi=-\log(\nu)$ satisfies
\begin{equation}\label{eqn:subharmonic}
0<2\tau^2\nu^4\leq\Delta\psi\leq -K\qquad\text{and}\qquad|\nabla\psi|\leq 2\tau\nu\sqrt{1-\nu^2}.
\end{equation}
In particular, $\psi$ is a strictly subharmonic function on $\Sigma_u$.
\end{lemma}

\begin{proof}
Let $\widetilde\kappa_1$ and $\widetilde\kappa_2$ be the eigenvalues of the shape operator of the dual entire $\tau$-graph $\widetilde\Sigma_v\subset\mathbb{L}^3$ at some point. Equation~\eqref{prop:K:eqn1} and the inequality $\widetilde\kappa_1^2+\widetilde\kappa_2^2\geq\frac{1}{2}(\widetilde\kappa_1+\widetilde\kappa_2)^2=2\tau^2$ give the estimate $\Delta\psi\geq 2\tau^2\nu^4$, whereas the estimate $\Delta\psi\leq -K$ follows from Lemma~\ref{lemma:psi}.(a) and the fact that $\widetilde K\leq 0$. Since $-2\tau\leq \widetilde\kappa_1,\widetilde\kappa_2\leq 0$ by Proposition~\ref{prop:CY}, the bound for the gradient is a consequence of Lemma~\ref{lemma:gradient}.
\end{proof}

Equation~\eqref{eqn:subharmonic} also says that the integral of $\Delta\psi$ over $\Sigma_u$ is bounded by the total curvature. If the latter is finite, then $\Sigma_u$ is conformally equivalent to $\mathbb{C}$, and we will apply standard estimates on the growth and the energy of $\psi$ by looking at it as a function on $\mathbb{C}$. In what follows, $\Delta_0$ and $\nabla_0$ will denote the usual Laplacian and gradient in $\C$, respectively, and $D_R(z_0)$ will stand for the open disk of center $z_0\in\C$ and radius $R>0$.

\begin{lemma}\label{lemma:growth}
Let $\psi\in\mathcal{C}^2(\C)$ be a nonnegative function such that $\Delta_0\psi\geq 0$ and
\[m=\int_\C\Delta_0\psi<\infty.\]
\begin{enumerate}[label=(\alph*)]
	\item There exists a constant $c>0$ such that $\psi(z)\leq c+\frac{2m}{\pi}\log(1+|z|)$ for all $z\in\C$.
	\item For any $R>0$, we have that $\int_{D_R(0)}|\nabla_0\psi|^2\leq 12\pi\sup_{D_{2R}(0)}\psi(z)^2$.
\end{enumerate}
\end{lemma}

\begin{proof}
Given $r>0$, let $M(r)=\frac{1}{2\pi}\int_0^{2\pi}\psi(re^{i\theta})\,\df\theta$ be the mean value of $\psi$ on the circle of radius $r$ centered at the origin. The divergence theorem yields
\[M'(r)=\frac{1}{2\pi}\int_0^{2\pi}\frac{\partial\psi}{\partial r}(re^{i\theta})\,\df\theta=\frac{1}{2\pi r}\int_{D_r(0)}\Delta_0\psi,\]
so that $0\leq M'(r)\leq\frac{m}{2\pi r}$. Therefore, the function $M$ is non-decreasing and satisfies $M(r)\leq M(1)+\frac{m}{2\pi}\log(r)$ for all $r\geq 1$. By the mean value inequality for subharmonic functions, we get
\[\psi(z)\leq\frac{1}{\pi|z|^2}\int_{D_{|z|}(z)}\psi\leq\frac{1}{\pi|z|^2}\int_{D_{2|z|}(0)}\psi=\frac{2}{|z|^2}\int_0^{2|z|}M(r)\,r\,\df r\leq 4M(2|z|),\]
for all $z\neq 0$, where we have used that $\psi\geq 0$ and $D_{|z|}(z)\subset D_{2|z|}(0)$. Consequently, if $|z|\geq\frac{1}{2}$, then $\psi(z)\leq 4M(1)+\frac{2m}{\pi}\log(2|z|)\leq 4M(1)+\frac{2m}{\pi}\log(2)+\frac{2m}{\pi}\log(1+|z|)$, and item (a) follows by taking $c$ as the maximum of $4M(1)+\frac{2m}{\pi}\log(2)$ and the maximum of $\psi$ on the closed disk of radius $\frac{1}{2}$ centered at the origin.

As for item (b), we will proceed as in a standard Caccioppoli type inequality, by considering the cut-off function 
\[\chi:\C\to[0,1],\qquad \chi(z)=\begin{cases}
1&\text{if }|z|\leq R,\\ 2-\frac{|z|}{R},&\text{if }R\leq|z|\leq 2R,\\ 0,&\text{if }2R\leq |z|.\end{cases}\]
Since $\chi^2\psi\nabla_0\psi$ has compact support,  integration by parts yields
\begin{equation}\label{lemma:growth:eqn2}
\begin{aligned}
\int_\C\chi^2|\nabla_0\psi|^2&=-\int_\C\chi^2\psi\,\Delta_0\psi-2\int_\C\chi\psi\,\langle\nabla_0\chi,\nabla_0\psi\rangle\\
&\leq\frac{1}{2}\int_\C\chi^2|\nabla_0\psi|^2+2\int_\C\psi^2|\nabla_0\chi|^2,
\end{aligned}
\end{equation}
where we have used that $\psi\,\Delta_0\psi\geq 0$ and the inequality $|\langle a,b\rangle|\leq\frac{1}{4}|a|^2+|b|^2$ applied to $a=\chi\nabla_0\psi$ and $b=\psi\nabla_0\chi$. Since $|\nabla_0\chi|=\frac{1}{R}$ in the annulus $R<|z|<2R$, whose area is $3\pi R^2$, and $\nabla_0\chi$ vanishes otherwise, we conclude from~\eqref{lemma:growth:eqn2} that
\[\int_{D_R(0)}\!\!|\nabla_0\psi|^2\leq\!\int_\C\!\chi^2|\nabla_0\psi|^2\leq 4\!\int_\C\!\psi^2|\nabla_0\chi|^2\leq 12\pi\sup_{D_{2R}(0)}\psi(z)^2.\qedhere\]
\end{proof}

\begin{proof}[Proof of Theorem~\ref{thm:FTC}]
Let us assume by contradiction that $\Sigma_u\subset\Nil(\tau)$ is an entire minimal graph with finite total curvature, and hence quadratic intrinsic area growth by~\cite[Thm.~3.4]{Man19}, so it is parabolic~\cite{CY1}. As $\Sigma_u$ is also simply connected, it is conformally equivalent to $\C$, i.e., there is a conformal parametrization $X:\C\to\Sigma_u$. If $\widetilde\Sigma_v\subset\mathbb{L}^3$ denotes the dual entire graph, then $\widetilde X=\Phi\circ X:\C\to\widetilde\Sigma_v$ is the dual conformal parametrization, where $\Phi$ is the conformal diffeomorphism given by~\eqref{eqn:phi}. As in Section~\ref{sec:gauss}, the metrics induced by $X$ and $\widetilde X$ on $\C$ will be written as $\mu^2|\df z|^2$ and $\widetilde\mu^2|\df z|^2$, respectively, where $\widetilde\mu=\nu\mu$. Now, all the functions defined on $\Sigma_u$ or $\widetilde\Sigma_v$ will be considered as functions on $\C$ through $X$ and $\widetilde X$. In particular, the Laplacians and the area elements of both surfaces are related to those of $\C$ by
\begin{equation}\label{thm:FTC:eqn1}
\Delta=\tfrac{1}{\mu^2}\Delta_0,\qquad\widetilde\Delta=\tfrac{1}{\widetilde\mu^2}\Delta_0,\qquad\df A=\mu^2\df A_0,\qquad\widetilde{\df A}=\widetilde\mu^2\df A_0.
\end{equation}

The function $\psi=-\log(\nu)$ is nonnegative and satisfies $0<\Delta\psi\leq|K|$ by Lemma~\ref{lemma:subharmonic}, so we deduce from~\eqref{thm:FTC:eqn1} that $\Delta_0\psi>0$ and
\begin{equation}\label{thm:TFC:eqn1.5}
m=\int_\C\Delta_0\psi\,\df A_0=\int_{\Sigma_u}\Delta\psi\,\df A\leq\int_{\Sigma_u}|K|\,\df A<\infty.
\end{equation}
Lemma~\ref{lemma:growth} applied to $\psi$ yields the existence of a constant $c>0$ such that
\begin{equation}\label{thm:FTC:eqn2}
\int_{D_R(0)}|\nabla_0\psi|^2\,\df A_0\leq12\pi\left(c+\tfrac{2m}{\pi}\log(1+2R)\right)^2,\quad\text{for all }R>0.
\end{equation}

The angle function $\widetilde\nu=-e^\psi$ lies in the kernel of the stability operator $\widetilde L$ given by~\eqref{eqn:jacobi}, so we get from~\eqref{thm:FTC:eqn1} that $\Delta_0\widetilde\nu=\widetilde\mu^2|\widetilde A|^2\widetilde\nu$, which can be written as
\begin{equation}\label{thm:FTC:eqn3}
\widetilde\mu^2|\widetilde A|^2=\frac{\Delta_0\widetilde\nu}{\widetilde\nu}=\Delta_0\psi+|\nabla_0\psi|^2.
\end{equation}
By integrating~\eqref{thm:FTC:eqn3} and taking~\eqref{thm:TFC:eqn1.5} and~\eqref{thm:FTC:eqn2} into account, we obtain
\begin{equation}\label{thm:FTC:eqn4}
\int_{D_R(0)}\widetilde\mu^2|\widetilde A|^2\,\df A_0\leq m+12\pi\left(c+\tfrac{2m}{\pi}\log(1+2R)\right)^2,\quad\text{for all }R>0.
\end{equation}
The Hopf differential $\widetilde Q\,\df z^2$ gives an entire function $\widetilde Q$ on $\mathbb{C}$ defined by~\eqref{eqn:Q}. By~\eqref{eqn:energy} and the arithmetic-geometric mean inequality, we get
\begin{equation}\label{thm:FTC:eqn5}
\widetilde\mu^2|\widetilde A|^2=2\tau^2\widetilde\mu^2+\frac{8|\widetilde Q|^2}{\widetilde\mu^2}\geq 8\tau|\widetilde Q|.
\end{equation}
Given $z_0\in\C$ and $R>0$, the mean value property for the holomorphic function $\widetilde Q$ on $D_R(z_0)\subseteq D_{R+|z_0|}(0)$, together with~\eqref{thm:FTC:eqn4} and~\eqref{thm:FTC:eqn5}, enables us to estimate
\[\begin{aligned}
|\widetilde Q(z_0)|&=\frac{1}{\pi R^2}\left|\int_{D_R(z_0)}\widetilde Q\,\df A_0\right|\leq\frac{1}{8\pi\tau R^2}\int_{D_{R+|z_0|}(0)}\widetilde\mu^2|\widetilde A|^2\,\df A_0\\
&\leq\frac{1}{8\pi\tau R^2}\left(m+12\pi\left(c+\tfrac{2m}{\pi}\log(1+2R+2|z_0|)\right)^2\right).
\end{aligned}\]
This last quantity tends to zero as $R\to\infty$, whence $\widetilde Q(z_0)=0$ for all $z_0\in\C$, i.e., $\widetilde\Sigma_v$ is a totally umbilical entire $\tau$-graph in $\mathbb{L}^3$, and hence a complete hyperbolic plane of constant curvature $-\tau^2$, namely the upper sheet of a hyperboloid of Lorentzian radius $\tau^{-1}$, which has hyperbolic conformal type. This is a contradiction.
\end{proof}

\begin{remark}
The invariant surfaces $\Sigma_{u_\theta}$ of Example~\ref{ex:FMP} show that the above estimates attain the equality in~\eqref{thm:FTC:eqn5}. For, it is easy to check that the change of variable $2\tau y=\sinh(2\tau t)$ defines a global conformal parameter $z=x+it$ on $\Sigma_{u_0}$, and hence on all the surfaces $\Sigma_{u_\theta}$ by Example~\ref{ex:boost}. In terms of this conformal parameter, we can set $\widetilde\mu=1$ and find that $\psi=\log\bigl(\cosh(\theta)\cosh(2\tau t)\bigr)$, which gives
\[\Delta_0\psi=\frac{4\tau^2}{\cosh^2(2\tau t)},\quad|\nabla_0\psi|^2=4\tau^2\tanh^2(2\tau t).\]
Thus $\psi$ grows linearly and $\widetilde\mu^2|\widetilde A|^2=4\tau^2$ is constant, and then the Hopf differential $\widetilde Q\,\df z^2=\frac{\tau}{2}\df z^2$ is also constant, so that~\eqref{thm:FTC:eqn5} becomes an equality.
\end{remark}

\section{Area, height and gradient estimates}\label{sec:growth}

In this section, we will prove Theorem~\ref{thm:height} and Corollaries~\ref{coro:area} and~\ref{coro:gradient}, as well as other related estimates. We will assume that $\tau>0$ in what follows. 

\subsection{Gradient estimates}
Our first goal will be to get a sharp gradient estimate for an entire spacelike $\tau$-graph in $\mathbb{L}^3$ under the assumption that it is bounded from below. The estimate will follow from analyzing the vertical sections of the graph, which turn out to be less curved than a certain hyperbola.

\begin{lemma}\label{lemma:sections}
Let $\widetilde\Sigma_v\subset\mathbb{L}^3$ be an entire spacelike $\tau$-graph. Given $p\in\R^2$ and a unit vector $w\in\R^2$, the function $\phi(s)=v(p+sw)$ satisfies $|\phi'(s)|<1$ and
\begin{equation}\label{eqn:sections}
0\leq\frac{\phi''(s)}{(1-\phi'(s)^2)^{3/2}}\leq 2\tau,\quad\text{for all }s\in\R.
\end{equation}
\end{lemma}

\begin{proof}
Consider the curve $\gamma(s)=(p+sw,\phi(s))$ contained in $\widetilde\Sigma_v$. Its derivative is given by $\gamma'=w_1\widetilde e_1+w_2\widetilde e_2$, where $w=(w_1,w_2)$ and the frame $\{\widetilde e_1,\widetilde e_2\}$ is defined in~\eqref{eqn:e1e2}. It follows from~\eqref{eqn:alpha-beta} and~\eqref{eqn:1ff} that
\begin{equation}\label{lemma:sections:eqn1}
\langle\gamma',\gamma'\rangle=1-(w_1\widetilde\alpha+w_2\widetilde\beta)^2=1-\phi'^2\geq1-\widetilde\alpha^2-\widetilde\beta^2=\widetilde\omega^2,
\end{equation}
where we have used that $\phi'=w_1v_x+w_2v_y$ and the Cauchy--Schwarz inequality. In particular, we get that $|\phi'|<1$. On the other hand, since $\phi''=w_1^2v_{xx}+2w_1w_2v_{xy}+w_2^2v_{yy}$, we can compute by means of~\eqref{eqn:2ff} the second fundamental form
\[\widetilde\sigma(\gamma',\gamma')=w_1^2\widetilde\sigma(\widetilde e_1,\widetilde e_1)+2w_1w_2\widetilde\sigma(\widetilde e_1,\widetilde e_2)+w_2^2\widetilde\sigma(\widetilde e_2,\widetilde e_2)=\frac{-\phi''}{\widetilde\omega}.\]
Hence~\eqref{eqn:CY-convex} can be written as $0\leq\phi''\leq 2\tau\,\widetilde\omega\,(1-\phi'^2)$, and the statement follows from the inequality $\widetilde\omega\leq(1-\phi'^2)^{1/2}$ given by~\eqref{lemma:sections:eqn1}.
\end{proof}

The curve $\gamma(s)=(p+sw,\phi(s))$ is the intersection of $\widetilde\Sigma_v$ and the vertical plane passing through $(p,v(p))$ containing $w$. This plane is isometric to the Minkowski plane $\mathbb{L}^2$, in which $\gamma$ is a spacelike curve that makes a Lorentzian angle $\vartheta=\arctanh(\phi')$ with the horizontal direction. The quotient in~\eqref{eqn:sections} is nothing but the curvature of $\gamma$ as a curve of $\mathbb{L}^2$, whence Lemma~\ref{lemma:sections} says that the vertical sections of $\widetilde\Sigma_v$ are convex curves with curvature at most $2\tau$, and this value is precisely the curvature of the hyperbola which is a vertical section orthogonal to the axis of the hyperbolic cylinder $v(x,y)=\frac{1}{2\tau}(1+4\tau^2y^2)^{1/2}$.

\begin{lemma}\label{lemma:comparison}
Let $\widetilde\Sigma_v\subset\mathbb{L}^3$ be an entire spacelike $\tau$-graph. If the function $v$ is bounded from below, then
\begin{equation}\label{eqn:comparison}
\widetilde\omega^{-1}\leq 1+2\tau\bigl(v-\inf\nolimits_{\R^2}v\bigr),
\end{equation}
where $\widetilde\omega=(1-v_x^2-v_y^2)^{1/2}$.
\end{lemma}

\begin{proof}
The inequality holds trivially at points where $\nabla v=0$ since the right-hand side is greater than or equal to $1$, so let us fix $p\in\R^2$ with $\nabla v(p)\neq 0$. Consider the function $\phi(s)=v(p+sw)$, where $w=-\frac{\nabla v(p)}{|\nabla v(p)|}$ is the direction of steepest descent of $v$ at $p$. Since $|\phi'|<1$, we can define the auxiliary function $\varsigma=\phi'(1-\phi'^2)^{-1/2}$ and Lemma~\ref{lemma:sections} yields the estimate $0\leq\varsigma'(s)\leq 2\tau$, which can be integrated to
\begin{equation}\label{prop:comparison:eqn1}
\varsigma(s)\leq\varsigma(0)+2\tau s=2\tau s-c,\quad\text{for all }s\geq 0,
\end{equation}
where
\begin{equation}\label{prop:comparison:eqn2}
c=-\varsigma(0)=\frac{-\phi'(0)}{\sqrt{1-\phi'(0)^2}}=\frac{|\nabla v(p)|}{\sqrt{1-|\nabla v(p)|^2}}>0.
\end{equation}
Since $\phi'=\varsigma(1+\varsigma^2)^{-1/2}$ and the function $t\mapsto t(1+t^2)^{-1/2}$ is increasing, we deduce from~\eqref{prop:comparison:eqn1} and~\eqref{prop:comparison:eqn2} (with $s=\frac{c}{2\tau}$) that
\[\begin{aligned}
v(p+\tfrac{c}{2\tau}w)-v(p)&=\int_0^{\tfrac{c}{2\tau}}\phi'(s)\,\df s\leq\int_0^{\tfrac{c}{2\tau}}\frac{(2\tau s-c)\,\df s}{\sqrt{1+(2\tau s-c)^2}}\\
&=\frac{1}{2\tau}\Bigl[\sqrt{1+(2\tau s-c)^2}\Bigr]_{s=0}^{s=\tfrac{c}{2\tau}}=\frac{1-\sqrt{1+c^2}}{2\tau}=\frac{1}{2\tau}\left(1-\frac{1}{\widetilde\omega(p)}\right).
\end{aligned}\]
As the left-hand side is at least $\inf_{\R^2}v-v(p)$, Equation~\eqref{eqn:comparison} follows readily.
\end{proof}


Lemma~\ref{lemma:comparison} does not give any information if $v$ is not bounded from below, e.g., for hyperbolic cylinders $\widetilde\Sigma_{v_\theta}$ given by~\eqref{eqn:FMP-dual} with $\theta\neq 0$. However, these examples suggest that the assumption can be dropped by applying a suitable boost, see Example~\ref{ex:boost}, which in turn gives a more geometric result (in the sense that the next statement is invariant under all isometries of $\mathbb{L}^3$ preserving the time orientation).

\begin{theorem}\label{thm:linear}
Let $\Sigma_u\subset\Nil(\tau)$ be an entire minimal graph and parametrize the dual entire $\tau$-graph $\widetilde\Sigma_v\subset\mathbb{L}^3$ (not necessarily conformally) as $\widetilde X:\R^2\to\mathbb{L}^3$ such that $\widetilde X(p)=(p,v(p))$, and with Gauss map $\widetilde N:\R^2\to\mathbb{H}^2$.
\begin{enumerate}[label=\emph{(\alph*)}]
	\item For all $p,q\in\R^2$, we have
	\begin{equation}\label{eqn:invariant}
	\langle\widetilde N(p),\widetilde N(q)\rangle\geq -1+2\tau\,\langle\widetilde X(p)-\widetilde X(q),\widetilde N(q)\rangle.
	\end{equation}
	\item Let $\delta:\R^2\to[0,\infty)$ be the function defined by $\delta=d_{\mathbb{H}^2}(\widetilde N,\widetilde E_3)$, which satisfies $\omega=\cosh(\delta)$. Then, for all $p,q\in\R^2$, we have
	\begin{equation}\label{eqn:linear}
	\omega(p)\leq e^{\delta(q)}\bigl(1+2\tau e^{\delta(q)}|p-q|\bigr),
	\end{equation}
	where $|p-q|$ denotes the Euclidean distance from $p$ to $q$.
\end{enumerate}
\end{theorem}

\begin{proof}
Let us fix $q\in\R^2$ and consider $\eta=\widetilde N(q)\in\mathbb{H}^2$, as well as the entire spacelike $\tau$-graph $\widetilde\Sigma_{v_\eta}=B_\eta(\widetilde\Sigma_v)$, where the boost $B_\eta$ is given by~\eqref{eqn:boost}, see Proposition~\ref{prop:deformation}. Note that $\widetilde\Sigma_{v_\eta}$ is parametrized by $B_\eta\circ\widetilde X$, but this parametrization transforms the base domain, and this transformation is given by the diffeomorphism $\varphi_\eta=\pi\circ B_\eta\circ\widetilde X:\R^2\to\R^2$. If we denote $\widetilde\omega_\eta=(1-|\nabla v_\eta|^2)^{1/2}$, then the height function and the angle function of $\widetilde\Sigma_{v_\eta}$ given by Proposition~\ref{prop:deformation}.(a) read
\begin{equation}\label{thm:linear:eqn1}
v_\eta\circ\varphi_\eta=-\langle\widetilde X,\eta\rangle\qquad\text{and}\qquad\frac{-1}{\widetilde\omega_\eta\circ\varphi_\eta}=\langle\widetilde N,\eta\rangle,
\end{equation}
respectively. In particular, we have that $\widetilde\omega_\eta(\varphi_\eta(q))=1$, which means that $\varphi_\eta(q)$ is a critical point of $v_\eta$. Since $v_\eta$ is a convex function by Proposition~\ref{prop:CY}, it attains its absolute minimum at $\varphi_\eta(q)$. Consequently, we can apply Lemma~\ref{lemma:comparison} to $\widetilde\Sigma_{v_\eta}$, so that~\eqref{eqn:comparison} evaluated at $\varphi_\eta(p)$ can be written by means of~\eqref{thm:linear:eqn1} as
\[-\langle\widetilde N(p),\eta\rangle\leq 1+2\tau\bigl(-\langle\widetilde X(p),\eta\rangle+\langle\widetilde X(q),\eta\rangle\bigr),\]
which is precisely~\eqref{eqn:invariant} and item (a) is proved. 

As for item (b), Proposition~\ref{prop:deformation} and Equation~\eqref{eqn:invariant} give
\begin{equation}\label{thm:linear:eqn2}
\omega(p)=-\langle\widetilde N(p),\widetilde E_3\rangle\leq -e^{\delta(q)}\langle\widetilde N(p),\widetilde N(q)\rangle\leq e^{\delta(q)}\bigl(1-2\tau\langle\widetilde X(p)-\widetilde X(q),\widetilde N(q)\rangle\bigr).
\end{equation}
By setting $w=\widetilde X(p)-\widetilde X(q)$ and $\eta=\widetilde N(q)$, we can use a similar reasoning as in~\eqref{prop:deformation:eqn1} to find that 
\begin{equation}\label{thm:linear:eqn3}
\bigl|\langle w,\eta\rangle+w_3\cosh(\delta(q))\bigr|\leq\sinh(\delta(q))(w_1^2+w_2^2)^{1/2}=\sinh(\delta(q))|p-q|.
\end{equation}
Moreover, $|w_3|=|v(p)-v(q)|\leq|p-q|$, which follows from the spacelike condition $|\nabla v|<1$. This, combined with~\eqref{thm:linear:eqn3}, yields
\begin{equation}\label{thm:linear:eqn4}
-\langle \widetilde X(p)-\widetilde X(q),\widetilde N(q)\rangle\leq(\sinh(\delta(q))+\cosh(\delta(q)))|p-q|=e^{\delta(q)}|p-q|,
\end{equation}
and the estimate~\eqref{eqn:linear} is a consequence of plugging~\eqref{thm:linear:eqn4} into~\eqref{thm:linear:eqn2}.
\end{proof}


\begin{remark}\label{rmk:linear}
Item (b) of Theorem~\ref{thm:linear} can be written in terms of the generalized gradient $Gu$ of the entire minimal graph $\Sigma_u\subset\Nil(\tau)$. Since $\omega^2=1+|Gu|^2$, we get that $|Gu|=\sinh(\delta)$ and $e^\delta=|Gu|+(1+|Gu|^2)^{1/2}$, so that~\eqref{eqn:linear} gives
\begin{equation}\label{eqn:linear-Gu}
|Gu(p)|<\omega(p)\leq e^{\delta(q)}+2\tau e^{2\delta(q)}|p-q|,\quad\text{for all }p,q\in\R^2.
\end{equation}
In particular, the estimate $|Gu|\leq C_1(1+r)$ in Theorem~\ref{thm:height} holds true for the constant $C_1=\max\{e^{\delta(0,0)},2\tau e^{2\delta(0,0)}\}$, which only depends on $\tau$ and $|Gu(0,0)|$. If the tangent plane of $\Sigma_u$ at $(q,u(q))$ is horizontal for some $q\in\R^2$, then $\delta(q)=0$ and~\eqref{eqn:linear} reduces to $\omega(p)\leq 1+2\tau|p-q|$. The minimal graph $u(x,y)=\tau xy$ shows that the constant $2\tau$ in this estimate cannot be improved, since it satisfies $\omega=(1+4\tau^2y^2)^{1/2}$ and has horizontal tangent plane at the origin.
\end{remark}


\begin{corollary}\label{coro:height}
Let $\Sigma_u\subset\Nil(\tau)$ be an entire minimal graph. Given $p=(p_1,p_2)\in\R^2$ and $q=(q_1,q_2)\in\R^2$, in the notation of Theorem~\ref{thm:linear}, we have that
\begin{equation}\label{eqn:height}
\bigl|u(p)-u(q)-\tau(q_1p_2-q_2p_1)\bigr|\leq e^{\delta(q)}|p-q|+\tau e^{2\delta(q)}|p-q|^2.
\end{equation}
\end{corollary}

\begin{proof}
From~\eqref{eqn:alpha-beta}, we get that $u_x=\alpha-\tau y$ and $u_y=\beta+\tau x$. Hence the derivative of $u$ along the segment $t\mapsto q+t(p-q)=(x(t),y(t))$ can be computed as
\[\begin{aligned}
\tfrac{\df}{\df t}u(q+t(p-q))&=(\alpha-\tau y)(p_1-q_1)+(\beta+\tau x)(p_2-q_2)\\
&=\alpha(p_1-q_1)+\beta(p_2-q_2)+\tau(q_1p_2-q_2p_1),
\end{aligned}\]
where $\alpha$ and $\beta$ are evaluated at $q+t(p-q)$. By integrating over $[0,1]$, we get
\begin{equation}\label{coro:height:eqn1}
u(p)-u(q)-\tau(q_1p_2-q_2p_1)=\int_0^1\bigl(\alpha(p_1-q_1)+\beta(p_2-q_2)\bigr)\df t.
\end{equation}
The Cauchy--Schwarz inequality and~\eqref{eqn:linear-Gu} imply that the absolute value of the integrand in~\eqref{coro:height:eqn1} is at most $|Gu|\,|p-q|\leq\bigl(e^{\delta(q)}+2\tau e^{2\delta(q)}t|p-q|\bigr)|p-q|$, whose integral over $[0,1]$ is the right-hand side of~\eqref{eqn:height}.
\end{proof}

\begin{remark}\label{rmk:height}
The left-hand side of~\eqref{eqn:height} is the height of the point $(p,u(p))$ after the left-translation that takes $(q,u(q))$ to the origin. Therefore, Corollary~\ref{coro:height} says that the height of an entire minimal graph with respect to any of its points grows at most quadratically. In particular, if $p=(x,y)$ and $q=(0,0)$, then~\eqref{eqn:height} gives
\[|u(x,y)|\leq|u(0,0)|+e^{\delta(0,0)}r+\tau e^{2\delta(0,0)}r^2,\quad\text{for all }p\in\R^2,\]
so the estimate $|u|\leq C_2(1+r^2)$ in Theorem~\ref{thm:height} holds true for a constant $C_2$ that depends on $\tau$, $|u(0,0)|$ and $|Gu(0,0)|$. This finishes the proof of Theorem~\ref{thm:height}.
\end{remark}

\begin{remark}\label{rmk:sharp}
The linear growth of the generalized gradient is sharp for any entire graph $\Sigma_u$ in $\Nil(\tau)$, not necessarily minimal. More specifically, for each $q\in\mathbb{R}^2$ there is some $p\in\mathbb{R}^2$ such that $|p-q|=r$ and $|Gu(p)|\geq\tau r$. 

To prove it, observe that $\operatorname{div}_{\mathbb{R}^2}(JGu)=\operatorname{div}_{\mathbb{R}^2}(\alpha\partial_y-\beta\partial_x)=\alpha_y-\beta_x=2\tau$, where $J$ is the $\frac{\pi}{2}$-rotation in $\mathbb{R}^2$. 
The divergence formula on $D_r(q)$ yields
\begin{equation}\label{eqn:heinz}
2\pi\tau r^2=\int_{\partial D_r(q)}\langle JGu,n\rangle\leq\int_{\partial D_r(q)}|Gu|\leq 2\pi r\max_{\partial D_r(q)}|Gu|,
\end{equation}
where $n$ is the unit conormal to $D_r(q)$ along its boundary. This Heinz type estimate gives no information if $\tau=0$, i.e., this property of the entire minimal graphs is a genuine feature of the non-zero bundle curvature.
\end{remark}

In~\cite{ManNel}, it was shown that $\nu^2$ is not integrable on an entire minimal graph, see also Remark~\ref{rmk:K}. Another consequence of~\eqref{eqn:linear} is the following improvement, which is sharp since $\nu^{3+\epsilon}$ is integrable on the umbrella of Example~\ref{ex:umbrella} for all $\epsilon>0$.

\begin{corollary}\label{coro:nu-not-L3}
Let $\Sigma_u\subset\Nil(\tau)$ be an entire minimal graph. Then $\nu^3\not\in L^1(\Sigma_u)$.
\end{corollary}

\begin{proof}
By a simple integration in polar coordinates, we have
\[\int_{\Sigma_u}\nu^3\,\df A=\int_{\R^2}\frac{\df x\,\df y}{\omega^2}\geq\int_0^{2\pi}\!\!\int_0^\infty\frac{r\,\df r\,\df\theta}{e^{2\delta(0,0)}\bigl(1+2\tau e^{\delta(0,0)}r\bigr)^2}=+\infty.\qedhere\]
\end{proof}

Finally, we will tackle the sharp gradient estimate for entire spacelike $\tau$-graphs in $\mathbb{L}^3$, which is just the counterpart to~\eqref{eqn:linear} by means of the duality.

\begin{corollary}\label{coro:gradient-explicit}
Let $\widetilde\Sigma_v\subset\mathbb{L}^3$ be an entire spacelike $\tau$-graph. Then
\begin{equation}\label{eqn:gradient}
|\nabla v(p)|^2\leq 1-\frac{1}{\lambda(q)^2\bigl(1+2\tau\lambda(q)\,|p-q|\bigr)^2},\quad\text{where }\lambda(q)^2=\frac{1+|\nabla v(q)|}{1-|\nabla v(q)|},
\end{equation}
for all $p,q\in\R^2$. In particular, the gradient estimate in Corollary~\ref{coro:gradient} holds true with the constant $C_3=\frac{1}{2\lambda(0,0)^2}\min\{1,\frac{1}{4\tau^2\lambda(0,0)^2}\}$.
\end{corollary}

\begin{proof} 
Let $\Sigma_u\subset\Nil(\tau)$ be the dual entire minimal graph. The inequality~\eqref{eqn:gradient} is a simple rewriting of~\eqref{eqn:linear} in the Lorentzian setting via duality. It is worth using that $|\nabla v|^2=1-\widetilde\omega^2=1-\omega^{-2}$ and $\lambda(q)=e^{\delta(q)}$, which follows from the fact that the distance $\delta(q)$ in Theorem~\ref{thm:linear} satisfies $\cosh(\delta(q))=\omega(q)$. Therefore,
\[\lambda(q)=e^{\delta(q)}=e^{\arctanh|\nabla v(q)|}=\sqrt{\frac{1+|\nabla v(q)|}{1-|\nabla v(q)|}}.\]

If $q=(0,0)$, then we can employ the inequality $(1+2\tau\lambda r)^2\leq 2(1+4\tau^2\lambda^2r^2)$, which follows from Cauchy--Schwarz inequality, in the denominator of~\eqref{eqn:gradient}. This easily leads to Corollary~\ref{coro:gradient} for the constant $C_3$ in the above statement. 
\end{proof}  

Note that, if the tangent plane of $\widetilde\Sigma_v$ at $(q,v(q))$ is horizontal for some $q\in\R^2$, then $\lambda=1$ and~\eqref{eqn:gradient} reduces to $1-|\nabla v(p)|^2\geq(1+2\tau|p-q|)^{-2}$. The hyperbolic cylinder $v(x,y)=\frac{1}{2\tau}(1+4\tau^2y^2)^{1/2}$, which has horizontal tangent plane at the origin and satisfies $1-|\nabla v|^2=(1+4\tau^2y^2)^{-1}$, shows that this estimate is asymptotically sharp as $|p-q|\to\infty$.

\subsection{Extrinsic area estimates} Given a surface $\Sigma$ immersed in $\Nil(\tau)$ and a point $p_0\in\Sigma$, denote by $B^\Sigma_R(p_0)$ the metric ball of $\Sigma$ with center $p_0$ and radius $R>0$, by $B_R(p_0)$ the metric ball of $\Nil(\tau)$ with center $p_0$ and radius $R$, and by $C_R(p_0)=\pi^{-1}(D_R(\pi(p_0)))$ the solid vertical cylinder of radius $R$ whose axis is the vertical line through $p_0$. Following~\cite{ManNel}, we will say that $\Sigma$ has \emph{exactly cubic} intrinsic (resp.\ \emph{extrinsic}, \emph{cylindrical}) \emph{area growth} if the area of $B^\Sigma_R(p_0)$ (resp.\ $\Sigma\cap B_R(p_0)$, $\Sigma\cap C_R(p_0)$) divided by $R^3$ lies between two positive constants for $R$ large enough, which does not depend upon the choice of $p_0$. These three notions of area are naturally ordered by the inclusions $B^\Sigma_R(p_0)\subset\Sigma\cap B_R(p_0)\subset\Sigma\cap C_R(p_0)$.

\begin{proof}[Proof of Corollary~\ref{coro:area}]
Let $\Sigma_u\subset\Nil(\tau)$ be an entire minimal graph. Since Theorem~\ref{thm:height} gives a constant $C_2>0$ such that $|u|\leq C_2(1+r^2)$, we deduce from~\cite[Cor.~3]{ManNel} that $\Sigma_u$ has exactly cubic extrinsic area growth. 

As for the cylindrical area growth, let $p_0\in\mathbb{R}^2$. On the one hand, it was proved in~\cite[Prop.~4]{ManNel} that $A(R)=\area(\Sigma_u\cap \pi^{-1}(D_R(p_0)))$ is at least the area of the corresponding piece of the umbrella centered at $p_0$, i.e., we can translate $p_0$ to the origin and use Example~\ref{ex:umbrella} to get the lower estimate
\begin{equation}\label{coro:area:eqn1}
A(R)\geq\int_0^{2\pi}\!\!\!\!\int_0^R\!\!r\sqrt{1+\tau^2r^2}\,\df r\,\df\theta=\tfrac{2\pi}{3\tau^2}\bigl((1+\tau^2R^2)^{3/2}-1\bigr)\geq\tfrac{2\pi\tau}{3}R^3.
\end{equation}
On the other hand, the area element of $\Sigma_u$ in the coordinates $(x,y)$ is $\omega\,\df x\,\df y$ in view of~\eqref{eqn:1ff}, so we can integrate~\eqref{eqn:linear} in polar coordinates centered at $p_0$ to get
\begin{equation}\label{coro:area:eqn2}
\begin{aligned}
A(R)&=\int_{D_R(p_0)}\!\!\!\!\omega\,\df x\,\df y\leq\int_0^{2\pi}\!\!\!\!\int_0^Re^{\delta(p_0)}\bigl(1+2\tau e^{\delta(p_0)}r\bigr)r\,\df r\,\df\theta\\
&=\pi e^{\delta(p_0)}R^2+\tfrac{4\pi\tau}{3}e^{2\delta(p_0)}R^3.
\end{aligned}
\end{equation}
It follows from~\eqref{coro:area:eqn1} and~\eqref{coro:area:eqn2} that $\Sigma_u$ has cubic cylindrical area growth.
\end{proof}


\begin{thebibliography}{88}

\bibitem{AR2}
U.\ Abresch, H.\ Rosenberg.
\newblock Generalized Hopf differentials. 
\newblock Proceedings of the XIII Escola de Geometria Diferencial.
\newblock \emph{Mat.\ Contemp.}, \textbf{28} (2005), no.\ 1, 1--28.

\bibitem{Cartier} S.\ Cartier.
\newblock Saddle towers in Heisenberg space.
\newblock \emph{Asian J.\ Math.}, \textbf{20} (2016), no.\ 4, 629--644.

\bibitem{CY1} S.\ Y.\ Cheng, S.\ T.\ Yau.
\newblock Differential equations on Riemannian manifolds and their geometric applications.
\newblock \emph{Comm.\ Pure Appl.\ Math.}, \textbf{28} (1975), no.\ 3, 333--354.

\bibitem{CY2} S.\ Y.\ Cheng, S.\ T.\ Yau.
\newblock Maximal space-like hypersurfaces in the Lorentz-Minkowski spaces.
\newblock \emph{Ann.\ of Math. (2)}, \textbf{104} (1976), no. 3, 407--419.

\bibitem{CT}  H.\ I.\ Choi, A.\ Treibergs.
\newblock Gauss maps of spacelike constant mean curvature hypersurfaces of Minkowski space.
\newblock \emph{J.\ Differential.\ Geom.}, \textbf{32} (1990), no.\ 3, 775--817.

\bibitem{CK} P.\ Collin, R.\ Krust.
\newblock Le probl\`{e}me de Dirichlet pour l'équation des surfaces minimales sur des domaines non bornés.
\newblock \emph{Bull.\ Soc.\ Math.\ France}, \textbf{119} (1991), no.\ 4, 443--462.

\bibitem{Dan} B.\ Daniel.
\newblock Isometric immersions into 3-dimensional homogeneous manifolds.
\newblock \emph{Comment.\ Math.\ Helv.}, \textbf{82} (2007), no.\ 1, 87--131.

\bibitem{Dan2} B.\ Daniel.
\newblock The Gauss map of minimal surfaces in the Heisenberg group.
\newblock \emph{Int.\ Math.\ Res.\ Not.}, \textbf{2011} (2011), no.\ 3, 674--695.

\bibitem{DFM} B.\ Daniel, I.\ Fernández, P.\ Mira.
\newblock The Gauss map of surfaces in $\widetilde{\mathrm{PSL}}_2(\mathbb{R})$.
\newblock \emph{Calc.\ Var.}, \textbf{52} (2015), no.\ 3--4, 507--528.

\bibitem{DH} B.\ Daniel, L.\ Hauswirth.
\newblock Half-space theorem, embedded minimal annuli and minimal graphs in the Heisenberg group.
\newblock \emph{Proc.\ Lond.\ Math.\ Soc. (3)}, \textbf{98} (2009), no.\ 2, 445--470.

\bibitem{DHM} B.\ Daniel, L.\ Hauswirth, P.\ Mira.
\newblock Lecture notes on homogeneous 3–manifolds.
\newblock 4th KIAS workshop on Differential Geometry, Seoul, 2009.

\bibitem{Esp} J.\ M.\ Espinar.
\newblock Finite index operators on surfaces.
\newblock \emph{J.\ Geom.\ Anal.}, \textbf{23} (2013), no.\ 1, 415--437.

\bibitem{FM} I.\ Fernández, P.\ Mira.
\newblock Holomorphic quadratic differentials and the Bernstein problem in Heisenberg space.
\newblock \emph{Trans.\ Amer.\ Math.\ Soc.}, \textbf{361} (2009), no.\ 11, 5737--5752.

\bibitem{FMP} C.\ Figueroa, F.\ Mercuri, R.\ H.\ L.\ Pedrosa.
\newblock Invariant surfaces of the Heisenberg groups.
\newblock \emph{Ann.\ Mat.\ Pura Appl. (4)}, \textbf{177} (1999), no.\ 1, 173--194.

\bibitem{Hartman} P.\ Hartman.
\newblock Geodesic parallel coordinates in the large.
\newblock \emph{Amer.\ J.\ Math.}, \textbf{86} (1964), no.\ 4, 705--727.

\bibitem{Lee} H.\ Lee.
\newblock Extensions of the duality between minimal surfaces and maximal surfaces.
\newblock\emph{Geom.\ Dedicata}, \textbf{151} (2011), no.\ 1, 373--386.

\bibitem{Li} P.\ Li.
\newblock Complete surfaces of at most quadratic area growth.
\newblock\emph{Comment.\ Math.\ Helv.}, \textbf{72} (1997), no.\ 1, 67--71.

\bibitem{Lopez} R.\ López.
\newblock Differential Geometry of curves and surfaces in Lorentz-Minkowski space.
\newblock \emph{Int.\ Electron.\ J.\ Geom.}, \textbf{7} (2014), no.\ 1, 44--107.

\bibitem{Man} J.\ M.\ Manzano.
\newblock On the classification of Killing submersions and their isometries.
\newblock \emph{Pac.\ J.\ Math.}, \textbf{270} (2014), no.\ 2, 367--392.

\bibitem{Man19} J.\ M.\ Manzano.
\newblock Dual quadratic differentials and entire minimal graphs in Heisenberg space.
\newblock \emph{Ann.\ Global Anal.\ Geom.}, \textbf{55} (2019), no.\ 2, 197--213.

\bibitem{ManNel} J.\ M.\ Manzano, B.\ Nelli.
\newblock Height and area estimates for constant mean curvature graphs in $\mathbb{E}(\kappa,\tau)$-spaces.
\newblock \emph{J.\ Geom.\ Anal.}, \textbf{27} (2017), no.\ 4, 3441--3473.

\bibitem{MPR} J.\ M.\ Manzano, J.\ Pérez, M.\ M.\ Rodríguez.
\newblock Parabolic stable surfaces with constant mean curvature.
\newblock\emph{Calc.\ Var.}, \textbf{42} (2011), no.\ 1--2, 137--152.

\bibitem{ManRod} J.\ M.\ Manzano, M.\ M.\ Rodríguez.
\newblock On complete constant mean curvature vertical multigraphs in $\mathbb{E}(\kappa,\tau)$.
\newblock \emph{J.\ Geom.\ Anal.}, \textbf{25} (2015), no.\ 1, 336--346.

\bibitem{MMPR} W.\ H.\ Meeks, P.\ Mira, J.\ Pérez, A.\ Ros.
\newblock Constant mean curvature spheres in homogeneous $3$-spheres.
\newblock \emph{J.\ Differential Geom.}, \textbf{120} (2022), no.\ 2, 307--343.

\bibitem{MeePer} W.\ H.\ Meeks, J.\ Pérez. 
\newblock Constant mean curvature surfaces in metric Lie groups. 
\newblock In: Geometric Analysis: Partial Differential Equations and Surfaces, \emph{Contemporary Mathematics} (AMS), \textbf{570} (2012), pp.~25--110.
\newblock

\bibitem{Milnor} T.\ K.\ Milnor.
\newblock Harmonic maps and classical surface theory in Minkowski 3-space.
\newblock \emph{Trans.\ Amer.\ Math.\ Soc.}, \textbf{280} (1983), no.\ 1, 161--185.

\bibitem{Treibergs} A.\ Treibergs.
\newblock Entire spacelike hypersurfaces of constant mean curvature in Minkowski space.
\newblock \emph{Invent.\ Math.}, \textbf{66} (1982), no.\ 1, 39--56.

\end{thebibliography}
\end{document}